\documentclass{article}
\usepackage{amscd}
\usepackage{amsmath}
\usepackage{amssymb}
\usepackage[mathscr]{euscript}
\usepackage{graphics}

\numberwithin{equation}{section}
\newtheorem{theorem}{Theorem}[section]

\newtheorem{proposition}[theorem]{Proposition}
\newtheorem{lemma}[theorem]{Lemma}
\newtheorem{example}[theorem]{Example}

\newcommand{\lam}{\lambda}
\newcommand{\Lam}{\Lambda}
\newcommand{\lz}{\lambda/\zeta}
\newcommand{\neta}{\nu/\eta}
\newcommand{\Sym}{{{\mathfrak S}}}
\newcommand{\Sn}{{{\mathfrak S}_n}}
\newcommand{\Hom}{{{\mathrm{Hom}}}}

\newcommand{\Ind}{{\mathrm{Ind}}}
\newcommand{\Res}{{\mathrm{Res}}}
\newcommand{\id}{{\mathrm{id}}}
\newcommand{\C}{{\Bbb C}}
\newcommand{\Z}{{\Bbb Z}}
\newcommand{\PP}{ {\mathcal P}}
\newcommand{\Pn}{ {\mathcal P}_n}
\newcommand{\Pic}{{\mathrm{Pic}}}
\newcommand{\PW}{{\mathrm{PW}}}
\newcommand{\PH}{{\mathrm{PH}}}
\newcommand{\uu}{{\bold u}}
\newcommand{\vv}{{\bold v}}
\newcommand{\xx}{{\bold x}}
\newcommand{\yy}{{\bold y}}
\newcommand{\zz}{{\bold z}}
\newcommand{\ww}{{\bold w}}
\newcommand{\uup}{{\bold u}^{\prime}}
\newcommand{\vvp}{{\bold v}^{\prime}}

\newcommand{\zzp}{{\bold z}^{\prime}}
\newcommand{\wwp}{{\bold w}^{\prime}}
\newcommand{\zetap}{\zeta^{\prime}}

\newcommand{\xip}{\xi^{\prime}}
\newcommand{\tp}{{\mathsf t}}
\newcommand{\leqnw}{\leq_{\nwarrow}}
\newcommand{\leqsw}{\leq_{\swarrow}}

\newcommand{\lsw}{<_{\swarrow}}
\newcommand{\Uplm}{\Upsilon_{\lambda \mu}}
\newcommand{\Psf}{{\mathsf \Pi}}
\newcommand{\LL}{{\mathsf \Lambda}}
\newcommand{\VV}{{\mathsf V}}
\newcommand{\IC}{\mathrm{IC}}
\newcommand{\ICC}{\mathrm{ICC}}
\newcommand{\bICC}{\overline{\mathrm{ICC} }}

\begin{document}



\title{A decomposition rule for certain tensor product representations of the symmetric groups }
\author{Takahiro Hayashi \quad}
%\email{hayashi@@math.nagoya-u.ac.jp}
\date{$\qquad$\\
%\date{%$\qquad$ 
Graduate School of Mathematics, 
Nagoya University, \\ %$\quad\;$
Chikusa-ku, Nagoya 464-8602, Japan\\
%\quad\\
%\quad 
email:\,\,hayashi@math.nagoya-u.ac.jp}
\maketitle

\begin{abstract}
In this paper, we give a combinatorial rule to calculate the decomposition 
of the tensor product (Kronecker product) of two irreducible complex representations of the   
symmetric group $\Sym_n$, when one of the representations corresponds to a hook $(n-m, 1^m)$.
\end{abstract}

%\begin{keyword}
%representations of symmetric groups 
%\sep 
%tensor product representations
%\sep
%Kronecker coefficient
%\MSC 05E10, 20C30



%\end{keyword}

%\end{frontmatter}


\section{Introduction}
\label{intro}

Let $G$ be a group such that each of its finite-dimensional complex representations 
is completely reducible.
One of the basic problems of the representation theory of $G$ is
to give a decomposition rule 
$$
 L(\lam) \otimes L(\nu) \cong \bigoplus_{\mu \in P} m^{\mu}_{\lam,\nu}   L(\mu)
$$
of the tensor product of two irreducible 
representations of $G$, where $\{ L(\lam)\, |\, \lam \in P\}$ denotes a set of complete 
representatives of irreducible   finite-dimensional representations of $G$.
When $G$ is the general linear group $GL(n,\C)$,  
the multiplicity $m^{\mu}_{\lam,\nu}  $ is known as the  famous 
Littlewood-Richardson coefficient $LR^{\mu}_{\lam,\nu} $,  
which equals to the number of  certain combinatorial objects 
(see e.g.  \cite{Fulton}). 
When $G$ is an arbitrary complex simple Lie group 
(or rather,  an arbitrary symmetrizable Kac-Moody algebra),  
inspired by Kashiwara's theory of crystals and 
works of Lakshmibai and Seshadri,  
Littelmann \cite{Littelmann} gave combinatorial objects
whose number equals to  $m^{\mu}_{\lam,\nu} $.

Let $G$ be the symmetric group $\Sym_n$ of $n$ letters.   
As usual, we use the set $\PP_n$  of partitions of $n$ as the set $P$ of labels of
irreducible representations of $G$.
In despite of the long research history (see e.g. Murnaghan  \cite{Murnaghan}
for an early work), 
much less is known about $m^{\mu}_{\lam,\nu} $ 
for $\Sym_n$-representations,
comparing with the Lie theoretic case.
Lascoux \cite{Lascoux},  Garsia and Remmel \cite{GarsiaRemmel},  
Remmel \cite{Remmel}, Remmel and Whitehead  \cite{RemmelWhitehead} 
and Rosas  \cite{Rosas} gave descriptions of 
$m^{\mu}_{\lam,\nu}$,  
when $\lam$ and $\nu$ are either two-row partitions or hook partitions. 
Recently, Ballantine and Orellana \cite{BO} gave a combinatorial rule for 
$m^{\mu}_{\lam,\nu}$ in the case where $\lam$ is a two-row partition $(n-p,p)$
and  $\nu$ is  not a partition inside the $2 (p-1) \times 2 (p-1)$ square.

In this paper, we give a combinatorial rule to calculate the  number
$m^{\mu}_{\lam,\nu} $ for $\Sym_n$-representations, 
when the partition $\nu$ is a hook $(n-m, 1^m)$.
More precisely,  we construct a set $\PH_m (\lambda, \mu)$ in a combinatorial manner,
which satisfies
$$
L(\lam) \otimes L (n-m, 1^m) \cong \bigoplus_{\mu \in \PP_n} |\PH_m (\lambda, \mu)| \cdot L(\mu)
$$
for each $\lambda, \mu \in \PP_n$ and $0 \leq m < n$.

Instead of dealing with the hook representation $L(n-m,1^m)$ directly, 
we consider a slightly bigger representation $ \Lam_m (\C^n)$, 
the $m$-th exterior power of the defining representation of $\Sym_n$.
By considering a certain 
permutation representation $\C \Omega_m$, 
we show that the multiplicity of
$L(\mu)$ in  $L(\lam) \otimes\Lambda_m (\C^n)$ is equal to  
$w_m:=\sum_{\zeta \in \PP_{n-m} } \sum_{\xi \in \PP_{m} } LR_{\zeta,\xi}^{\lambda}$ 
$\cdot LR_{\zeta,\xi^{\tp}}^{\mu}$,
where $\xi^{\tp}$ denotes the transpose of  $\xi$.
Although we could not find the result in the literature,
we suspect that it was known for experts, 
since the techniques to be used to prove is rather standard. 

To give a decomposition of $L(\lam) \otimes L (n-m, 1^m)$, we 
use a set $\PW_m (\lambda, \mu)$ such that $|\PW_m (\lambda, \mu)| = w_m$ and that
each of its elements is a  Zelevinsky's {\it picture} \cite{Zelevinsky}.
A picture is a bijective map between skew Young diagrams, which satisfies certain order-theoretic  conditions, and it is identified with a tableau on a skew Young diagram satisfying  some conditions.
Using a variant of  Zelevinsky's insertion algorithm for pictures \cite{Zelevinsky}, 
we construct the set  $\PH_m (\lambda, \mu)$ as a subset of  $\PW_m (\lambda, \mu)$.

\section{Preliminaries}
\label{Preliminaries}
\subsection{Partitions and diagrams}
\label{Partitions and diagrams}

Let  $\lambda = (\lam_1, \lam_2, \ldots, \lam_l)$ be a sequence of  integers.
We say that $\lam$ is a {\it partition} of $n$ 
if  $\lam_1 \geq  \ldots \geq  \lam_l > 0$ and $| \lam | := \sum_i \lam_i = n$. 
For a partition
$\lambda = (\lam_1, \ldots, \lam_l)$, 
we define its {\it length} $l(\lam)$   by
$l (\lam ) = l$.
Also, we define its {\it diagram} $D(\lam)$ and its $i$-{\it th row}  
$D_i (\lambda) \quad (1 \leq i \leq l)$ by
$D(\lambda) = \amalg_i D_i (\lam)$ and
$D_i (\lambda) = \{ (i, 1), (i,2), \ldots, (i,\lambda_i) \}$, respectively.  
We denote by  $\Pn$  the set of partitions of $n$.
For a partition $\lam$ and $i \geq 1$, we define $\lam_i \in \Z_{\geq 0}$
by $\lam = (\lam_1, \lam_2, \ldots, \lam_{l (\lam)})$
and $ \lam_i = 0\,\,\, (i > l (\lam))$.

For two partitions  $\lam \in \Pn$ and $\zeta \in \PP_m $, we write $\zeta \subseteq \lam$ if 
$D(\zeta) \subseteq D(\lam)$.
When $\zeta \subseteq \lam$, we denote the pair $(\lam,\zeta)$ by
$\lam/\zeta$ and call it a {\it skew partition} of $(n,m)$.
For a skew partition $\lam/\zeta$, we define its {\it size}, its {\it length}, its {\it diagram}
and its $i$-{\it th row}
by $| \lam/\zeta | := | \lam | - | \zeta |$, $l( \lam/\zeta ) := l( \lam )$, 
$D(\lam/\zeta): = D(\lam) \setminus D(\zeta)$ and 
$D_i (\lam/\zeta): = D_i (\lam) \setminus D_i (\zeta)$,
respectively.

For a partition $\lam$, we define its {\it conjugate} $\lam^{\tp}$
to be a unique partition satisfying 
$D(\lam^{\tp}) = D(\lam)^{\tp} := \{ \uu^{\tp}\,|\, \uu  \in  D(\lam) \}$, 
where  
$(i,j)^{\tp} := (j,i)$. For a skew partition $\lam/\zeta$, we set
$(\lam/\zeta)^{\tp} := \lam^{\tp} /\zeta^{\tp} $.

\subsection{Representations}
\label{Representations}

Let $G$ be a finite group. 
Let $M$ be a finite-dimensional $\C G$-module 
and let $L(\lam)$ be a simple $\C G$-module labeled by $\lam$.
We denote by $[M]$, $[M: L(\lam)]$ and $M_{\lam}$ the character of $M$, 
the multiplicity of $L(\lam)$ in $M$ and the homogeneous component 
of $M$ corresponding to $L(\lam)$, respectively. 

%************Lemma\\
\begin{lemma}
\label{FrobRecip}
Let $G$ be a finite group and let $H$ be its subgroup.
For each $\C G$-module $M$ and 
$\C H$-module $N$, we have the  isomorphisms
 $M \otimes \Ind^G_H N \cong \Ind^G_H (\Res^G_H M \otimes N);$
$x \otimes (\sigma \otimes_{\C H}  y) \mapsto \sigma \otimes_{\C H}  (\sigma^{-1} x \otimes y)\,\,$
$(x \in M, \sigma \in G, y \in N)$ of $\C G$-modules.
Also, we have the isomorphism  
$\Hom_{\C H} (\Res^G_H M, N) \cong \Hom_{\C G} (M, \Ind^G_H N);$   
$g \mapsto f;$ $f(x) = \frac{1}{|H|} \sum_{\sigma \in G}  \sigma  \otimes_{\C H} g (\sigma^{-1} x)$
$(x \in M, g \in \Hom_{\C H} (\Res^G_H M, N))$ of vector spaces.
Here $ \Ind^G_H$ and  $\Res^G_H $ denote the induction functor and the restriction functor, 
respectively.
\end{lemma}

See e.g. \cite{CurtisReiner} 
Corollary 10.20 and Proposition 10.21 for a proof.\\

Let $\Sn$ be the symmetric group of $n$ letters.
For each $\lambda \in \Pn$, we denote the corresponding simple $\C \Sn$-module
by $L(\lambda)$ (see e. g.  \cite{Fulton} or  \cite{GoodmanWallach}).
It is well-known that
\begin{equation}
\label{L^*}
 L(\lam)^* \cong L(\lam),
\quad
\C_{sgn} \otimes L(\lam) \cong L(\lam^{\tp}), 
\end{equation}
where $L(\lam)^*$ denotes the dual module of $L(\lam)$ and 
$\C_{sgn}$ denotes the module corresponding to the sign representation.
For each $\sigma \in \Sn$, we denote by
$\sigma_{\lam} \in GL(L(\lam))$ the 
image of $\sigma$ via the representation corresponding to $L(\lam) $.

Let $\Omega = \Omega_{1,n}$ be the set $\{1,2,\ldots, n \}$ equipped with the natural 
action of $\Sn$. Then the {\it defining module}  $\C \Omega \cong {\C}^n$ has the following decomposition:
$$
\C \Omega \cong L(n) \oplus L(n-1,1).  
$$
Moreover, we have
$$
 \Lambda_m (L(n-1,1)) \cong L(n-m,1^m)$$
for each $0 \leq m \leq n-1$
(see e.g. \cite{GoodmanWallach} page 391),  where $\Lam_m$ stands for the  $m$-th exterior power
and $(n-m,1^m)$ stands for the {\it $m$-th hook}
$$
(n-m,1^m) =(n-m,\overbrace{1,1,\ldots,1}^{m {\text{times}}}).
$$
Since $L(n)$ is isomorphic to the unit module (trivial module) $\C_{unit}$, we also have
\begin{equation}
\label{LamCOmega}
 \Lambda_m (\C \Omega) \cong 
\begin{cases}
L(n) & (m= 0)\\
L(n - m + 1,1^{m-1}) \oplus L(n-m,1^m) & (0 < m < n)\\
L(1^n) & (m = n). \\
\end{cases}
\end{equation}
Let $\#_m\!: \Sym_m \to \Sn$ be the group homomorphism which sends each
transposition 
$(i, i+1)$ to $(n - m + i, n - m + i + 1)$.
As usual, we identify the subgroup 
$\langle \rho \cdot \#_m (\tau)\,|\, \rho \in \Sym_{n-m}, \tau \in \Sym_m \rangle$ 
of $\Sn$ with 
$\Sym_{n-m} \times\Sym_m$.   
The {\it Littlewood-Richardson coefficient} $LR^{\lam}_{\zeta,\xi} $ 
is defined
by  the following decomposition rule:  
\begin{equation}
\label{DefLR}
[\mathrm{Res}^{\Sn}_{\Sym_{n - m} \times \Sym_m} L(\lam)]
=
\sum_{ \zeta \in \PP_{n - m}}  \sum_{\xi \in \PP_{m}} 
LR^{\lam}_{\zeta,\xi} [L(\zeta) \boxtimes L(\xi) ], 
\end{equation}
where $\boxtimes$ stands for the outer tensor product.

Let $M$ and $N$ be $\C \Sn$-modules and let $g\!: M \to N$ be a $\C \Sn$-module map.
For each $\lam, \mu \in \Pn$, we define a vector space $\Uplm (M)$ and a linear map  
$
\Uplm (g)\!: \Uplm (M)  \to \Uplm (N) 
$
by
$$
 \Uplm (M) = \Hom_{\C \Sn}
 (L(\mu), L(\lambda) \otimes M)  
$$
and $\Uplm (g)  (f) = (\id_{L(\lam)} \otimes g) \circ f \quad (f \in \Uplm (M))$.
Then $\Uplm$ gives an exact functor from the category of 
finite-dimensional $\C\Sn$-modules to the category of finite-dimensional 
$\C$-vector spaces.  We note that
\begin{equation}
\label{Mult=DimUp}
[L(\lam) \otimes M]=  \sum_{\mu \in \Pn} \dim (\Uplm (M)) \cdot [L(\mu)].
\end{equation}

\section{The module $\C \Omega_m$}
\label{ThemoduleCOmegam}

For each $0 < m \leq n$, we consider the subset $\Omega_m = \Omega_{m,n}$ 
of the $m$-fold Cartesian product $\Omega^m$ defined by
$$
\Omega_m 
:=
\bigl\{ (i_1,\ldots, i_m) \in \Omega^m \bigm|\, i_1,\ldots,i_m\,\, 
\text{are pairwise distinct} \bigr\}.
$$
We define an action 
$(\sigma, I) \mapsto \sigma_V I \quad (\sigma \in \Sn, I \in \Omega_m)$ 
of $\Sn$ on $\Omega_m$
by
\begin{equation}
  \sigma_V (i_1, \ldots, i_m) = (\sigma (i_1), \ldots, \sigma (i_m)) 
\end{equation}
and call it the {\it vertical action}.
Also, we define the  {\it horizontal action}
$(\tau, I) \mapsto \tau_H I \quad (\tau \in \Sym_m, I \in \Omega_m)$ 
by
$$
  \tau_H (i_1, \ldots, i_m) = (i_{\tau^{-1} (1)}, \ldots, i_{\tau^{-1} (m)}) .
$$
Since these actions commute with each other, the linear span
$\C \Omega_m$ becomes an $\Sn\times\Sym_m$-module. 

It is easy to see that the vertical action of $\Sn$ on $\Omega_m$
is transitive and that the stabilizer of $I_m: = (n -m + 1, n - m + 2, \ldots, n)$
is identified with $\Sym_{n-m}$. 
Hence $\C \Omega_m$ is isomorphic to the
induced module $\Ind_{\Sym_{n-m}}^{\Sn} L(n-m)$.

For each  $f \in \Uplm ({\Bbb C} \Omega_m )$ and
$I \in \Omega_m$, we define $f_I\!:\,L(\mu) \to L(\lambda)$ by 

$$
f(v) 
= 
\sum_{I}
f_I (v) \otimes I.
$$

The following proposition follows immediately from 
Lemma 
\ref{FrobRecip}.

\begin{proposition}
For each $\lam, \mu \in \Pn$, there exists a linear isomorphism
\begin{equation}
 \label{Heartsuit}
 \heartsuit_m\!:
 \Hom_{{\Bbb C}{\mathfrak S}_{n - m}} (L(\mu), L(\lambda))  
\cong
 \Uplm ( {\Bbb C} \Omega_m ) 
\end{equation}
defined by
\begin{multline}
\qquad\quad  \heartsuit_m  (h) (v) = \frac{1}{ (n - m)!} \sum_{\sigma \in \Sn}  \sigma h(\sigma^{-1} v) \otimes \sigma_V I_m\\
(h \in \Hom_{{\Bbb C}{\mathfrak S}_{n - m}} (L(\mu), L(\lambda)) , v \in L(\mu) ).\quad \quad 
\end{multline}
The inverse of $\heartsuit_m$ is given by 
$$
\heartsuit_m  ^{-1} (f) = f_{I_m}
\quad 
(f \in \Uplm ({\Bbb C} \Omega_m )).
$$
\end{proposition}

Since $\Uplm$ is a functor,  
$\Uplm( {\Bbb C} \Omega_m )$ becomes  an $\Sym_{m}$-module via 
$\tau \mapsto \Uplm (\tau_H)$.
On the other hand, the left-hand side of \eqref{Heartsuit} also becomes 
an $\Sym_m$-module via  
$\tau_{\lam \mu} (h) := \#_m (\tau)_{\lam} \circ h \circ \#_m (\tau)_{\mu}^{-1} \,\,$
$(h \in  \Hom_{{\Bbb C}{\mathfrak S}_{n - m}} (L(\mu), L(\lambda)) ,\, \tau \in \Sym_m)$.  
In fact, we have the following:
\begin{lemma}
The map $\heartsuit_m$ commutes with the actions of $\Sym_m$.
In particular, we have the following linear isomorphism of homogeneous
components:
\begin{equation}
 \label{IsoHomComp1}
 \Hom_{{\Bbb C}{\mathfrak S}_{n - m}} (L(\mu), L(\lambda)) _{(1^m)} 
\cong
 \Uplm ( {\Bbb C} \Omega_m )_{(1^m)} .
\end{equation}
%

%
%\begin{equation}
%\label{CDHeart}
%\begin{CD}
%\Hom_{{\Bbb C}{\mathfrak S}_{n - m}} (L(\mu), L(\lambda))   
%@ > \heartsuit_m  >> \Uplm ( {\Bbb C} \Omega_m ) \\
% @ V \tau_{\lam \mu} VV @ VV  \Uplm (\tau_H) V \\
% \Hom_{{\Bbb C}{\mathfrak S}_{n - m}} (L(\mu), L(\lambda))  
% @ > \heartsuit_m  >> \Uplm ( {\Bbb C} \Omega_m ) \\
%\end{CD}
%\end{equation}
%
\end{lemma}

\noindent
{\it Proof.}
To prove the assertion,  it suffices to show
\begin{align}
\label{Heart&Sm}
 \heartsuit_m^{-1} (\Uplm (\tau_H) (f) )
= &
\tau_{\lam, \mu} (\heartsuit_m^{-1}  (f))
\end{align}
for each transposition 
$\tau = (p, p+1) \in \Sym_m$
and $f \in \Uplm (\C \Omega_m)$.
For each $I = (i_1, i_2, \ldots, i_m) \in \Omega_m$, 
we have
\begin{equation}
\label{tauHI}
\tau_H^{-1} I
= 
(i_1, \ldots, i_{p-1}, i_{p+1}, i_{p}, i_{p+2},\ldots, i_m)
=
(i_p,i_{p+1})_V I.
\end{equation}
On the other hand, 
comparing  
$$
 \sigma f(v) 
 = \sum_J \sigma_{\lam} f_J (v) \otimes \sigma_V J
 = \sum_I \sigma_{\lam} f_{ \sigma_V^{-1}  I} (v) \otimes I
$$
with
$$
 \sigma f(v) =  f(\sigma_{\mu} v) 
 = \sum_I f_I (\sigma_{\mu} v) \otimes I, 
$$
for each $\sigma \in \Sn$  and $v \in L(\mu)$, we get 
\begin{equation}
\label{f_sV(I)}
 f_{\sigma_V^{-1} I} (v) = \sigma_{\lam}^{-1} f_I  (\sigma_{\mu} v). 
\end{equation}
In particular, 
\begin{equation*}
\label{f_tH(I)}
 f_{\tau_H^{-1} I} (v) 
 = 
 (i_p, i_{p+1})_{\lam}
  f_I \left((i_p, i_{p+1})_{\mu}^{-1} v \right)
\end{equation*}
by \eqref{tauHI}. Hence,  we have
\begin{align*}
(\Uplm (\tau_H) (f)) (v)
&=
\sum_I f_{\tau_H^{-1} I} (v) \otimes I\\
&= 
\sum_{I = (i_1,\ldots,i_m)}  (i_p, i_{p+1})_{\lam}
  f_I \left((i_p, i_{p+1})_{\mu}^{-1} v \right) \otimes I.
\end{align*}
Since $(i_p,i_{p+1}) = \#_m (p, p + 1)$ for $I = I_m$, we get
\eqref{Heart&Sm}.
$\hfill \square$

Let $(\C \Omega_m)_{- \boxtimes (1^m)}$ be the homogeneous 
component of $\C \Omega_m$ with respect to
the horizontal action of $\Sym_m$, which corresponds to $L(1^m)$.
It is easy to see that the correspondence   
$(i_1) \wedge\cdots \wedge(i_m) \mapsto$ $ \sum_{\tau} \mathrm{sgn} (\tau)$ 
$ (i_{\tau^{-1} (1)}, \ldots, i_{\tau^{-1} (m)})$ 
gives an isomorphism 
\begin{equation}
\label{LamCong}
 \Lam_m (\C \Omega)
 \cong 
(\C \Omega_m)_{- \boxtimes (1^m)}
\end{equation}
of $\C \Sn$-modules.

\begin{lemma}
We have the following linear isomorphisms:
\begin{equation}
\label{IsoHomComp2}
 \Uplm (\C \Omega_m)_{(1^m)}
 \cong 
 \Uplm ((\C \Omega_m)_{- \boxtimes (1^m)}),
\end{equation}
\begin{equation}
\label{IsoHomComp3}
\left( L(\mu)^* \otimes L(\lam)\right)_{(n-m)\boxtimes (1^m)}  
\cong 
\Hom_{\C \Sym_{n-m}}  (L(\mu), L(\lam))_{(1^m)},
\end{equation}
where the left-hand side of the second isomorphism
denotes the homogeneous component of
$\Res^{\Sn}_{\Sym_{n-m} \times \Sym_m} (L(\mu)^* \otimes L(\lam))$, 
which corresponds to %the one-dimensional module
%simple $\C(\Sym_{n-m} \times \Sym_m)$-module
$L(n-m) \boxtimes L(1^m)$.
\end{lemma}
\medskip
\par\noindent
{\it Proof}. 
Define an action of $\Sn \times \Sym_m$ on 
$H:= \Hom_{\C} (L(\mu), L(\lam) \otimes \C \Omega_m)$
by 
$(\sigma,\tau)f =
( \sigma_\lam \otimes (\sigma_V \tau_H)) \circ f \circ\sigma_{\mu}^{-1} $.
Then the left-hand side of \eqref{IsoHomComp2} is
identified with the homogeneous component
$H_{(n)\boxtimes (1^m)}$, 
which corresponds to 
$L(n) \boxtimes L(1^m)$.
Let $\iota\!:(\C \Omega_m)_{- \boxtimes (1^m)} \to \C \Omega_m$
be the natural injection. 
To prove  \eqref{IsoHomComp2}, it is enough to show that the subspace
$\mathrm{Im} (\Uplm (\iota) )$ of $H$ coincides with $H_{(n)\boxtimes (1^m)}$.
Let $f$ be an element of $H_{(n)\boxtimes (1^m)}$.
Since $f(v) \in (L(\lam) \otimes \C \Omega_m)_{- \boxtimes (1^m)}
=  L(\lam) \otimes (\C \Omega_m)_{- \boxtimes (1^m)}$
for each $v \in L(\mu)$, $f$ gives a linear map
$f^{\prime}\!:L(\mu) \to L(\lam) \otimes (\C \Omega_m)_{- \boxtimes (1^m)}$.
Since both $\iota$ and $f$ commute with the actions of $\Sn$, 
$f^{\prime}$ also commutes with the actions of $\Sn$. This proves   
$\mathrm{Im} (\Uplm (\iota)) \supseteq H_{(n)\boxtimes (1^m)}$.
Since the other  inclusion "$\subseteq$" is obvious, we have completed the proof of
\eqref{IsoHomComp2}. 
The isomorphism  \eqref{IsoHomComp3} is rather obvious, 
since both of the spaces in  \eqref{IsoHomComp3} are
naturally identified with 
$\Hom_{\C}(L(\mu),L(\lam))_{(n-m) \boxtimes (1^m)}$.
$\hfill \square$

\begin{lemma}
We have
\begin{equation}
\label{DimHomComp}
\dim \left( L(\mu)^* \otimes L(\lam)\right)_{(n-m)\boxtimes (1^m)}  
=
\sum_{\zeta \in \PP_{n-m} } \sum_{\xi \in \PP_{m} } 
LR_{\zeta,\xi}^{\lambda} LR_{\zeta,\xi^{\tp}}^{\mu}.
\end{equation}
\end{lemma}
\medskip
\par\noindent
{\it Proof}. 
By  \eqref{L^*} and \eqref{DefLR}, we have
%
%\begin{multline}
\begin{gather}
[\mathrm{Res}^{\Sn}_{\Sym_{n - m} \times \Sym_m} 
\left( L(\mu)^* \otimes L(\lam) \right)\, :\, L(n-m) \boxtimes L(1^m)  ]\nonumber \\
=
\sum_{ \zeta, \zetap \in \PP_{n - m}}  \sum_{\xi,\xip \in \PP_{m}} 
LR^{\lam}_{\zeta,\xi} LR^{\mu}_{\zetap,\xip}
\qquad\qquad\qquad\qquad\nonumber \\
\label{ResLmuLLam}
 \qquad\qquad
\times[L(\zetap) \otimes  L(\zeta)\,:\, L(n-m) ]
 [L(\xip)  \otimes L(\xi)\,:\,  L(1^m) ].
\end{gather}
%\end{multline}
%
On the other hand, by \eqref{L^*}, we have
$$
 [L(\zetap) \otimes  L(\zeta)\,:\, L(n-m) ]
=
\dim \Hom_{\C \Sym_{n-m}}( L(\zetap)^*, L(\zeta))
= \delta_{\zeta \zetap}
$$
and
$$
 [L(\xip) \otimes  L(\xi)\,:\, L(1^m) ]
=
\dim \Hom_{\C \Sym_{m}}( L(1^m) \otimes L(\xip)^*, L(\xi))
= \delta_{\xi^{\tp} \xip}.
$$
Hence the left-hand side of \eqref{ResLmuLLam} agrees with 
the right-hand side of \eqref{DimHomComp}.
Since $L(n - m) \boxtimes L(1^m)$ is one-dimensional, this proves \eqref{DimHomComp}.
$\hfill \square$

\begin{theorem}
\label{L(lam)Lam}
 $\mathrm{(1)}$ The multiplicity $w_m$ of $L(\mu)$ in the tensor product 
 module $L(\lam) \otimes\Lambda_m (\C\Omega)$ is equal to  
\begin{equation}
\sum_{\zeta \in \PP_{n-m} } \sum_{\xi \in \PP_{m} } LR_{\zeta,\xi}^{\lambda} LR_{\zeta,\xi^{\tp}}^{\mu}.
\end{equation}
$\mathrm{(2)}$ The multiplicity $h_m$ of $L(\mu)$ in 
%the tensor product module 
 $L(\lam) \otimes L(n - m, 1^m)$ is equal to  
$\sum_{i = 0}^m (-1)^{m - i} w_i$.
\end{theorem}

\medskip
\par\noindent
{\it Proof}. 
The number $w_m$ is equal to $\dim \Uplm (\C \Omega_m)$ by \eqref{Mult=DimUp}.
Then, applying \eqref{LamCong}, 
\eqref{IsoHomComp2}, \eqref{IsoHomComp1}, \eqref{IsoHomComp3}
and \eqref{DimHomComp} in this order, Part (1) follows.
By \eqref{Mult=DimUp}, 
 $h_m$ is equal to $\dim \Uplm (L(n-m, 1^m))$.
Hence, we have
$w_0 = h_0$, $w_m= h_{m-1} + h_m\,\, (0 < m <n)$ and $w_n = h_{n-1}$
by \eqref{LamCOmega}. 
Solving these equations, we get Part (2). $\hfill \square$

\section{A crystallized exact sequence}
\label{A crystallized exact sequence}

Let  $\lam/\zeta$ and $\nu/\eta$ be skew partitions of  $(n,n - m)$.
In \cite{Zelevinsky}, Zelevinsky has constructed a finite set $\Pic (\nu/\eta, \lam/\zeta)$
of {\it pictures}, 
which satisfies 
$$
|\Pic (\nu/\eta, \lam/\zeta)| = \sum_{\xi \in \PP_m}   
LR_{\zeta, \xi}^{\lam} 
\cdot   LR_{\eta,\xi }^{\nu} .
$$ 
By Theorem \ref{L(lam)Lam} (1) and the well-known symmetry 
$LR_{\zeta^{\tp},\omega^{\tp}}^{\mu^{\tp}} = LR_{\zeta,\omega}^{\mu}$,
the multiplicity of  
$L(\mu)$ in  $L(\lam) \otimes\Lambda_m (\C\Omega)$ is
equal to the number of elements of the set
\begin{equation}
\label{DefPEm}
\PW_m (\lambda, \mu) := \amalg_{\zeta  \in \PP_{n - m} } \PW (\lambda, \mu; \zeta),
\end{equation}
where
\begin{equation}
\label{DefPElmz}
 \PW (\lambda, \mu; \zeta) := \Pic((\mu/\zeta)^{\tp}, \lam/\zeta)
\end{equation}
if $\zeta \subseteq \lam, \mu$, and $\PW (\lambda, \mu; \zeta) = \emptyset$
if otherwise.
In this section, we will construct a subset 
$\PH_m (\lambda, \mu) \subseteq \PW_m (\lambda, \mu) $
 whose number of elements  is equal to
the 
multiplicity of  
$L(\mu)$ in  $L(\lam) \otimes L(n - m,1^m)$.
Also, we will give a bijection 
$E\!: \PH_m (\lambda, \mu) \cong 
\PW_{m+1} (\lambda, \mu)\setminus\PH_{m+1} (\lambda, \mu) $.
For this purpose, we will use a backward analogue of Zelevinsky's insertion-deletion
algorithm for pictures.
As a first step, we introduce a backward analogue of Sagan-Stanley's 
insertion-deletion algorithm for  tableaux on skew partitions    
(cf. \cite{SaganStanley}).

\bigskip\par\noindent
{\it Remark.}
By \eqref{LamCOmega}, there exists an exact sequence
$$
\Bbb{E}: 0 \to \Lam_0 (\C\Omega) \to \Lam_1 (\C\Omega)
\to \cdots \to \Lam_n (\C\Omega) \to 0
$$
of $\C \Sn$-modules.
Since $\Upsilon_{\lam\mu}$ is exact, there also exists an exact sequence
$\Upsilon_{\lam\mu} (\Bbb{E})$ of vector spaces.
The above-mentioned construction may be viewed as a combinatorial counterpart
of $\Upsilon_{\lam\mu} (\Bbb{E})$.
\subsection{Tableaux}
\label{Tableaux}

For a partition $\lambda  \in \Pn$,
we say that a point $\uu  \in \Z^2_{>0}$ is a {\it corner} of 
$\lam$ if  $\uu \in D(\lam)$ and $D(\lam) \setminus \{ \uu \} = D(\xi)$ for 
some $\xi =: \lam \setminus \uu \in \PP_{n-1}$.
Also, 
we say that a point $\vv \in \Z^2_{>0}$ is a {\it cocorner} of 
$\lam$ if  $\vv \not\in D(\lam)$ and $D(\lam) \amalg \{ \vv \} = D(\omega)$ for 
some $\omega =: \lam \amalg \vv \in \PP_{n+1}$.
For a skew partition $\lambda/\zeta$,
we say that a point $\ww \in D(\lambda/\zeta)$ is an {\it inner corner} of 
$\lz$ if  $\ww$ is a cocorner of $\zeta$.
We denote by  $\IC (\lambda/\zeta)$ the set of inner corners of 
$\lambda/\zeta$.
We say that a point $\zz  \in \Z^2_{>0}$ is an {\it inner cocorner} of 
$\lambda/\zeta$ if  $\zz$ is a corner of $\zeta$.
We say that $\zz \in \Z^2$ is an {\it extreme  cocorner} of $\lambda/\zeta$
if either $\zz = (l(\lambda/\zeta),0)$ and 
$(l(\lambda/\zeta),1) \in D(\lambda/\zeta)$, or
$\zz = (0, \lambda_1)$ and 
$(1,\lambda_1) \in D(\lambda/\zeta)$.
Let $\mathrm{ICC} (\lambda/\zeta)$ be the set of inner cocorners of $\lambda/\zeta$
and let $\overline{\mathrm{ICC} }(\lambda/\zeta)$ be the set of inner cocorners and 
extreme cocorners of $\lambda/\zeta$.
For example, when 
$\lambda/\zeta = (5,5,4,2,1)/(3,3,2,1)$, 
we have $\IC (\lambda/\zeta) = \{ \ww_1, \ww_2, \ww_3,\ww_4 \}$, 
$\ICC (\lambda/\zeta) = \{ \zz_1, \zz_2,\zz_3\}$ and
$\bICC (\lambda/\zeta) = \{ \zz_i\,|\, 0 \leq i \leq 4 \}$,
where $\ww_1 = (5,1)$, $\ww_2, \ldots,\zz_4$ are as in the first figure
of \eqref{ExampleIC} below.   
When 
$\lambda/\zeta = (5,5,4,2,1)/(5,3,2,1,1)$, 
we have $\IC (\lambda/\zeta) = \{ \ww_1, \ww_2, \ww_3 \}$
and 
$\ICC (\lambda/\zeta) = \{ \zz_0, \zz_1, \zz_2,\zz_3\} = \bICC (\lambda/\zeta) $, 
where $\ww_1$, $\ww_2, \ldots,\zz_3$ are as in the second figure
of \eqref{ExampleIC}.

\begin{equation}
\label{ExampleIC}
\setlength{\unitlength}{1mm}
\begin{picture}(40,27)
 \put(10,25){\line(1,0){25}}
 \put(20,15){\line(1,0){5}}
 \put(30,15){\line(1,0){5}}
 \put(15,10){\line(1,0){10}}
 \put(25,10){\line(1,0){5}}
 \put(10,5){\line(1,0){10}}
 \put(10,0){\line(1,0){5}}
 \put(10,0){\line(0,1){25}}
 \put(15,5){\line(0,1){5}}
 \put(20,10){\line(0,1){5}}
 \put(25,15){\line(0,1){10}}
 \put(15,0){\line(0,1){5}}
 \put(20,5){\line(0,1){5}}
 \put(30,10){\line(0,1){5}}
 \put(35,15){\line(0,1){10}}
 %
 %\put(10,15){\framebox(4.9,4.9){$\vv_1$}}
% \put(15,15){\framebox(4.9,4.9){}} 
 \put(5,0){\makebox(4.9,4.9){$\zz_0$}}
 \put(10,0){\makebox(4.9,4.9){$\ww_1$}}
 \put(10,5){\makebox(4.9,4.9){$\zz_1$}}
 \put(15,5){\makebox(4.9,4.9){$\ww_2$}}
 \put(15,10){\makebox(4.9,4.9){$\zz_2$}}
 \put(20,10){\makebox(4.9,4.9){$\ww_3$}}
 \put(20,15){\makebox(4.9,4.9){$\zz_3$}}
 \put(25,20){\makebox(4.9,4.9){$\ww_4$}}
 \put(30,25){\makebox(4.9,4.9){$\zz_4$}}
\end{picture}
\qquad
\setlength{\unitlength}{1mm}
\begin{picture}(40,27)
 \put(10,25){\line(1,0){25}}
 \put(25,20){\line(1,0){10}}
 \put(20,15){\line(1,0){5}}
 \put(30,15){\line(1,0){5}}
 \put(15,10){\line(1,0){10}}
 \put(25,10){\line(1,0){5}}
 \put(15,5){\line(1,0){5}}
 \put(10,0){\line(1,0){5}}
 \put(10,0){\line(0,1){25}}
 \put(15,5){\line(0,1){5}}
 \put(20,10){\line(0,1){5}}
 \put(25,15){\line(0,1){5}}%
 \put(15,0){\line(0,1){5}}
 \put(20,5){\line(0,1){5}}
 \put(30,10){\line(0,1){5}}
 \put(35,15){\line(0,1){10}}
 \put(10,0){\makebox(4.9,4.9){$\zz_0$}}
 \put(15,5){\makebox(4.9,4.9){$\ww_1$}}
 \put(15,10){\makebox(4.9,4.9){$\zz_1$}}
 \put(20,10){\makebox(4.9,4.9){$\ww_2$}}
 \put(20,15){\makebox(4.9,4.9){$\zz_2$}}
 \put(25,15){\makebox(4.9,4.9){$\ww_3$}}
 \put(30,20){\makebox(4.9,4.9){$\zz_3$}}
\end{picture}
\end{equation}

\medskip\medskip
Let
$\leqnw$ and 
$\leqsw$ be partial orders on $\Z^2$ defined by
\begin{align}
 (i,j) \leqnw (k,l) &
 \Longleftrightarrow i\leq k\, \text{and}\, j \leq l,\\
 (i,j) \leqsw (k,l) &
 \Longleftrightarrow k\leq i\, \text{and}\, j \leq l.
\end{align}
We note that
\begin{align}
\label{Tenchi&Leqsw}
 \xx \leqsw \yy &
 \Longleftrightarrow \yy^{\tp}  \leqsw \xx^{\tp}
\end{align}
for each $\xx, \yy \in \Z^2$.
Also, we note that
$(\mathrm{IC}  (\lambda/\zeta)\amalg \bICC (\lambda/\zeta), \leqsw)$ 
is a totally ordered set
whose order is expressed as
\begin{align}
\label{z0w1zkwk}
 \zz_0 \lsw \ww_1 \lsw \zz_1 \lsw \ww_2 
 \lsw \ldots \lsw \ww_k \lsw \zz_k,
\end{align}
where 
$\IC (\lambda/\zeta) = \{ \ww_1, \ldots, \ww_k \}$
and
$\bICC (\lambda/\zeta) = \{ \zz_0, \ldots, \zz_k \}$.

Let $\nu/\eta$ be a skew partition and let $T\!: D( \nu/\eta) \to \Z_{>0}$ be an injective map.
We say that $T$ is a  {\it partial tableau} on  $\nu/\eta$ if it is a map of ordered sets from 
$(D( \nu/\eta), \leqnw)$ to $(\Z_{>0}, \leq)$.
%Let $a > 0$ be an integer. 
%We define a point $\uu_l \in \{ l \} \times \Z$ as follows: 
%If there exists no element $\uu \in D_l (\nu/\eta)$ satisfying 
%$T(\uu) < a$, then put $\uu_l = (l, \eta_l) \not\in D_l (\nu/\eta)$.
%If not, then define $\uu_{l}$
%to be the right-most point in $D_l (\nu/\eta)$ satisfying 
%$T(\uu_{l}) < a$. 
%When $\uu_l \in D_l (\nu/\eta)$, we also define 
%$\uu_{l-1} \in \{ l-1 \} \times \Z$ to be $(l-1, \eta_{l-1})$  if 
%there exists no element $\uu \in D_{l-1} (\nu/\eta)$ satisfying $T(\uu) < T(\uu_{l-1}) $ 
Let  $T$ be a partial tableau on a skew partition $\nu/\eta$ of length
$l = l (\nu/\eta)$. 
Let $a > 0$ be an integer such that $a \not\in T(D( \nu/\eta))$. 
In order to construct the backward row insertion of $a$ into $T$, 
we define the {\it  bumping route} $\xx_l, \xx_{l-1}, \ldots, \xx_r$
of  $(T, a)$ by the following lemma:

\begin{lemma}
\label{RowInstBumpRoute}
There exist a unique integer $l \geq r \geq 0$ and a unique sequence  
$\xx_l, \xx_{l-1}, \ldots, \xx_r \in \Z^2$ satisfying 
the following three conditions:

\begin{quote}
{\rm (1)} $\xx_i \,\, (i  > r)$ is the right-most point of  $D_i (\nu/\eta)$
which satisfies $T(\xx_i) <  a_{i} $, 
where  $a_j =T (\xx_{j+1})$ for  $r \leq j<l$ and $a_l = a$.\\
{\rm (2)}  If $r > 0$, then $a_r\leq T(\xx)$ for every $\xx \in D_r (\nu/\eta)$.\\
{\rm (3)} Either $r > 0$ and $\xx_r = (r, \eta_{r})$,  or $r = 0$
and $\xx_0 = (0,\nu_1)$. 
\end{quote}
\end{lemma}

The lemma is proved by determining $\xx_l, \xx_{l-1}, \ldots, \xx_r$ inductively. 
We call $\uu:= \xx_r \in \bICC (\neta)$ 
the {\it bumping destination} of $(T,a)$ and 
write $\displaystyle{\uu\, \frac{\,\,\,T\,\,\,}{}\, a}$.
We say that $a$ is an {\it addable integer} of $T$ 
if $\uu \in \mathrm{ICC} (\neta)$.
In this case, we define the {\it backward} {\it  row insertion} $E_{a} T$ 
to be the partial tableau on 
$\nu/(\eta \setminus \uu)$ determined by 
$(E_{a} T) (\xx_i) = a_{i} \,\,\,\, (l\geq  i \geq r)$ 
and 
$(E_{a} T)(\yy) = T(\yy) \,\,\,\,  (\yy \ne \xx_l, \xx_{l-1}, \ldots, \xx_r)$.

As well as the usual row insertion, 
the backward row insertion has an inverse operation. 
Let $\vv$ be an inner corner of $\nu/\eta$.
Then there exists at most one pair $(T^{\prime}, b)$ such that 
$T = E_b T^{\prime}$ and that $\vv$ is the bumping destination of $(T^{\prime}, b)$.
When there exists such a pair, 
we say that $\vv$ is a  {\it removable corner} of $T$ and we write
$T^{\prime} = F_{\vv} T$ and $\displaystyle{\vv\, \frac{}{\,\,\,T\,\,\,}\, b}$. 
More explicitly, $F_{\vv} T$ is given as follows:

%\bigskip\bigskip
%Next,  let $\uu \in \{r \} \times \Z$ be an inner corner of $\nu/\eta$.
%We say that $\uu$ is a {\it removable corner} of $T$ 
%if these exists a (necessarily unique)
%sequence $\uu_r = \uu, \uu_{r + 1} , \ldots, \uu_l$  of points of $D(\nu/\eta)$ 
%such that
%$\uu_i$ is the left-most point in $D_i (\nu/\eta)$ satisfying 
%$T(\uu_{i-1}) < T(\uu_i)$ for each $r < i \leq l$.  
%In this case, we call  $\uu_r, \uu_{r + 1} , \ldots, \uu_l$
%the {\it bumping route} of $( \uu,T)$.
%Also, we write $\displaystyle{\uu\, \frac{}{\,\,\,T\,\,\,}\, T(\uu_l)}$
%and define the {\it row deletion} $F_{\uu} T$ to be the partial tableau on 
%$\nu/(\eta \amalg  \uu)$ determined by 
%$(F_{\uu} T) (\uu_i) = T(\uu_{i - 1}) \,\,\,\, (r < i < l)$ 
%and 
%$(F_{\uu} T)(\vv) = T(\vv) \,\,\,\,  (\vv \ne \uu_{r+1},  \ldots, \uu_l)$.

\begin{lemma}
\label{RowDelBumpRoute}
An inner corner  $\vv = (s,j)$ of $\neta$ is a  removable corner of $T$ 
if and only if these exists a (necessarily unique)
sequence $\yy_s = \vv, \yy_{s + 1} , \ldots, \yy_l$  of points of $D(\nu/\eta)$ 
such that
$\yy_i$ is the left-most point in $D_i (\nu/\eta)$ satisfying 
$T(\yy_{i-1}) < T(\yy_i)$ for each $s < i \leq l$.  
In this case, we have 
$(F_{\vv} T) (\yy_i) = T(\yy_{i-1})\,\,(s < i \leq l)$, 
$(F_{\vv} T) (\xx) = T(\xx)\,\,(\xx \ne \yy_{s+1},\ldots, \yy_l)$ and 
$\displaystyle{\vv\, \frac{}{\,\,\,T\,\,\,}\, T(\yy_l)}$.
\end{lemma}

%\medskip%\medskip
%\par\noindent
%{\it Remark.} The above algorithm also works for semistandard tableaux.\\
%without any change.

\begin{example}
\label{ExampleT}
{\rm 
Define a partial tableau $T$ on $(5,5,4,3)/(4,3,2)$
as follows:
\begin{equation}
\label{ExampleTEq}
\setlength{\unitlength}{0.8mm}
\begin{picture}(15,10)
 \put(0,15){\makebox(4.9,4.9){$T=$}}
\end{picture}
\setlength{\unitlength}{0.8mm}
\begin{picture}(15,35)
 \put(0,5){\line(0,1){20}}
 \put(0,25){\line(1,0){25}}
 \put(0,5){\framebox(4.9,4.9){$1$}}
 \put(5,5){\framebox(4.9,4.9){$5$}}
 \put(10,5){\framebox(4.9,4.9){$10$}}
 \put(10,10){\framebox(4.9,4.9){$3$}}
 \put(15,10){\framebox(4.9,4.9){$11$}}
 \put(15,15){\framebox(4.9,4.9){$8$}}
 \put(20,15){\framebox(4.9,4.9){$9$}}
 \put(20,20){\framebox(4.9,4.9){$6$}}
\end{picture}
\setlength{\unitlength}{0.8mm}
\begin{picture}(40,20)
 \put(10,15){\makebox(4.9,4.9){$.$}}
\end{picture}
\end{equation}
Then, the bumping route of $(T,7)$ is 
$(4,2)$,  $(3,3)$, $(2,3)$.
Hence $\displaystyle{(2,3)\, \frac{\,\,\,T\,\,\,}{}\, 7}$,
and
$7$ is an addable integer of $T$.
The tableau $E_7 T$ is given by
the first equality of \eqref{ExampleET} below.
\begin{equation}
\label{ExampleET}
\setlength{\unitlength}{0.8mm}
\begin{picture}(10,10)
 \put(0,15){\makebox(4.9,4.9){$ E_{7}T=$}}
\end{picture}
\setlength{\unitlength}{0.8mm}
\begin{picture}(20,35)
 \put(0,5){\line(0,1){20}}
 \put(0,25){\line(1,0){25}}
 \put(0,5){\framebox(4.9,4.9){$1$}}
 \put(5,5){\framebox(4.9,4.9){$7$}}
 \put(10,5){\framebox(4.9,4.9){$10$}}
 \put(10,10){\framebox(4.9,4.9){$5$}}
 \put(15,10){\framebox(4.9,4.9){$11$}}
 \put(10,15){\framebox(4.9,4.9){$3$}}
 \put(15,15){\framebox(4.9,4.9){$8$}}
 \put(20,15){\framebox(4.9,4.9){$9$}}
 \put(20,20){\framebox(4.9,4.9){$6$}}
 \put(25,15){\makebox(4.9,4.9){$,$}}
\end{picture}
\setlength{\unitlength}{0.8mm}
\begin{picture}(40,20)
 \put(15,15){\makebox(4.9,4.9){$\qquad\qquad\qquad\qquad F_{(3,3)}T=\qquad\qquad$}}
\end{picture}
\setlength{\unitlength}{0.8mm}
\begin{picture}(15,35)
 \put(0,5){\line(0,1){20}}
 \put(0,25){\line(1,0){25}}
 \put(0,5){\framebox(4.9,4.9){$1$}}
 \put(5,5){\framebox(4.9,4.9){$3$}}
 \put(10,5){\framebox(4.9,4.9){$10$}}
 \put(15,10){\framebox(4.9,4.9){$11$}}
 \put(15,15){\framebox(4.9,4.9){$8$}}
 \put(20,15){\framebox(4.9,4.9){$9$}}
 \put(20,20){\framebox(4.9,4.9){$6$}}
\end{picture}
\end{equation}
On the other hand, $\vv = (3,3)$ is a removable corner of $T$
and $F_{(3,3)} T$ is given by the second equality of  \eqref{ExampleET}.
Moreover, we have
$\displaystyle{(3,3)\, \frac{}{\,\,\,T\,\,\,}\, 5}$.

}
\end{example}

\medskip\medskip
The following "bumping lemmas" are quite essential for our discussions.

\begin{lemma}
\label{BumpLemTab}
For each partial tableau $T$ on $\nu/\eta$, we have the following:\\
 $\mathrm{(1)}$  Let $a \leq a^{\prime}$ be positive integers such that $a,a^{\prime} \not\in T(D(\neta))$.
Define $\uu, \uup \in \bICC(\nu/\eta)$ by 
$\displaystyle{\uu\, \frac{\,\,\,T\,\,\,}{}\, a}$
and 
$\displaystyle{\uup\, \frac{\,\,\,T\,\,\,}{}\, a^{\prime}}$.
Then, we have $\uu \leqsw \uup$. \\
%
%
%$\mathrm{(2)}$ Let $\vv$ and $\vvp$ be elements of $\IC (\nu/\eta)$
%such that 
%$\vv \leq_{\swarrow} \vvp$ and that 
%$\vvp$ is a removable corner of $T$.
%Then $\vv$ is also a removable corner of $T$.
%Moreover,  we have $b \leq b^{\prime}$, where
%$\displaystyle{\vv\, \frac{}{\,\,\,T\,\,\,}\, b}$ and $\displaystyle{\vvp\, \frac{}{\,\,\,T\,\,\,}\, b^{\prime}}$. \\
%
$\mathrm{(2)}$ Let $\vv$ and $\vvp$ be removable corners of $T$ such that
$\vv \leqsw \vvp$.
Define integers $b$ and $b^{\prime}$ by 
$\displaystyle{\vv\, \frac{}{\,\,\,T\,\,\,}\, b}$ and $\displaystyle{\vvp\, \frac{}{\,\,\,T\,\,\,}\, b^{\prime}}$. 
Then, we have $b \leq b^{\prime}$.\\
$\mathrm{(3)}$  Let $a$ be an addable integer of $T$ and let  $a^{\prime}$ be another positive integer
such that $a^{\prime} \not \in T(D(\neta)) \amalg \{ a \}$.
Define $\uu\in \ICC(\nu/\eta)$ and 
$\uup \in  \bICC(\nu/(\eta\setminus \uu))$ by 
$\displaystyle{\uu\, \frac{\,\,\,T\,\,\,}{}\, a}$
and 
$\displaystyle{\uup\, \frac{\,\,\,E_a T\,\,\,}{}\, a^{\prime}}$, respectively.
Then we have
$a < a^{\prime} \Rightarrow \uu \lsw \uup$ and $a^{\prime} < a \Rightarrow \uup \lsw \uu$.\\
%
%
%$\mathrm{(4)}$  Let $\vv$ be a removable corner of $T$ with
%$\displaystyle{\vv\, \frac{}{\,\,\,T\,\,\,}\, b}$
%and let  $\vvp$ be an inner corner of 
%$\nu/(\eta \amalg \vv)$.
%If $\vvp \lsw \vv$, then $\vvp$ is also removable and satisfies $b^{\prime} < b$, where 
%$\displaystyle{\vvp\, \frac{}{\,\,\,F_{\vv} T\,\,\,}\, b^{\prime}}$.
%If $\vv \lsw \vvp$ and $\vvp$ is also removable, then  $b < b^{\prime}$, where 
%$\displaystyle{\vvp\, \frac{}{\,\,\,F_{\vv} T\,\,\,}\, b^{\prime}}$.
%
$\mathrm{(4)}$
Let $\vv$ be a removable corner of $T$ and 
let $\vvp$ be a removable corner of $F_{\vv}T$.
Define integers $b$ and $b^{\prime}$ by 
$\displaystyle{\vv\, \frac{}{\,\,\,T\,\,\,}\, b}$ and 
$\displaystyle{\vvp\, \frac{}{\,\,\,F_{\vv} T\,\,\,}\, b^{\prime}}$, respectively. 
Then, we have 
$\vv \lsw \vvp \Rightarrow b < b^{\prime}$ and $\vvp \lsw \vv \Rightarrow b^{\prime} < b$.\\
%
%\par\noindent
%$\mathrm{(5)}$
%Let $\uu$ be an inner cocorner of $\nu/\eta$
%and let $\vv$  be an inner corner of $\nu/\eta$.
%Suppose that there exists an addable integer $a$ of $T$ such that
%$\displaystyle{\uu\, \frac{\,\,\,T\,\,\,}{}\, a}$.
%If $\vv \lsw \uu$, then 
% $\vv$ is a removable corner of $T$ and satisfies  $b < a$, where
%$\displaystyle{\vv\, \frac{}{\,\,\,T\,\,\,}\, b}$.
%If $\uu \lsw \vv$ and 
% $\vv$ is a removable corner of $T$ with $\displaystyle{\vv\, \frac{}{\,\,\,T\,\,\,}\, b}$,  
%then $a < b$.
%
%
$\mathrm{(5)}$
Let $a \not\in T(D(\neta))$ be a positive integer and
define $\uu \in \bICC (\nu/\eta)$ by 
$\displaystyle{\uu\, \frac{\,\,\,T\,\,\,}{}\, a}$.
Let $\vv$ be a removable corner of $T$ such that $\uu \lsw \vv$. 
Then, we have $a < b$, where $\displaystyle{\vv\, \frac{}{\,\,\,T\,\,\,}\, b}$.\\
$\mathrm{(6)}$
Let $a$ and $\uu$ be as in Part (5).
If a cocorner $\vv = (s,j)$ of $\eta$ satisfies $\vv \lsw \uu$, then
$\vv$ is a removable corner of $T$.
Moreover, we have $b < a$, where  
$\displaystyle{\vv\, \frac{}{\,\,\,T\,\,\,}\, b}$.
\end{lemma}

\noindent
{\it Proof.} 
All of these can be proved in a  standard manner  (cf. \cite{Fulton} page 9).
Here we will prove Part (6) using Lemma \ref{RowDelBumpRoute}.
Let $\xx_l, \xx_{l-1}, \ldots , \xx_r = \uu$ be the 
bumping route of $(T,a)$. 
Since $s > r$, we have $D_s (\nu/\eta) \supseteq \{\xx_s \} \ne \emptyset$.
Hence $\vv$ is an inner corner of $\nu/\eta$. 
Suppose there exist an integer $s \leq t < l$ and 
a  sequence  $\yy_s: = \vv, \yy_{s + 1} , \ldots, \yy_t$, 
such that  each $\yy_i\,\, (s < i \leq t)$ satisfies the condition in Lemma  \ref{RowDelBumpRoute} 
and that $\yy_t$ is left of or equal to  $\xx_t$.
Since $T(\yy_t) \leq T(\xx_t) < T(\xx_{t+1})$, 
we have $S:= \{ \yy \in D_{t + 1} (\neta)\,|\,T(\yy_t) < T(\yy) \} \ne \emptyset$.
Let  $\yy_{t+1}$ be the left-most element of $S$.  
Then, it is obvious that $\yy_{t+1}$ satisfies the condition in 
Lemma  \ref{RowDelBumpRoute} for $i = t + 1$ and that
$\yy_{t+1}$ is left of or equal to  $\xx_{t+1}$.
By induction, we see that there exists a sequence  $\yy_s, \yy_{s + 1} , \ldots, \yy_l$
satisfying the condition of Lemma \ref{RowDelBumpRoute}
and that $\yy_l$  is left of or equal to  $\xx_l$.
This proves the first assertion as well as the second assertion, since 
$b = T( \yy_l ) \leq T(\xx_l ) < a$.
$\hfill \square$

\subsection{Pictures}
\label{Pictures}

Next, we recall the definition of the picture of Zelevinsky, following 
Fomin and Greene \cite{FominGreene} (or rather, Leeuwen \cite{Leeuwen}).

%Let $\chi$ be a finite non-empty subset of $\Z^2$. We say that $\chi$ is a {\it skew shape}
%if it is $\leqnw$-convex, that is, 
%if $\uu \leqnw \vv \leqnw \zz$ and $\uu, \zz \in \chi$, then $\vv \in \chi$.  

%Let $\chi$ and $\psi$ be skew shapes and let $p\!: \chi \to \psi$ be a bijection.
%We say that $p$ is a {\it picture} from $\chi$ onto $\psi$ if 
%$\uu, \vv \in \chi$, $\uu \leqsw \vv$ $p(\uu) \leqnw p(\vv)$
%$p \in \Pic (\chi, \psi)$ 

%$|\Pic (\chi, \psi)| = \sum_{\lam \in \Pn} |\Pic (\chi,\lam)|  \cdot  |\Pic (\lam, \psi)|  $  

%$|\Pic(\mu, \lam/\zeta)| = LR^{\lam}_{\mu,\zeta}$

Let $\nu/\eta$ and $\lam/\zeta$ be skew partitions  of $(n,n-m)$.
We say that a bijection $\Psf\!: D(\nu/\eta) \to D(\lam/\zeta)$
is a {\it picture} from $\nu/\eta$ onto $\lam/\zeta$ if both $\Psf\!: (D(\nu/\eta),\leqnw) \to (D(\lam/\zeta),\leqsw)$ 
and $\Psf^{-1}\!: (D(\lam/\zeta),\leqnw) \to (D(\nu/\eta),\leqsw)$ 
are maps of ordered sets.
We denote by $\Pic (\nu/\eta, \lam/\zeta)$  
the set of pictures  from $\nu/\eta$ onto $\lam/\zeta$.

Let
$R\!:(D(\lam/\zeta),\leqsw) \to \Z_{>0}$ be an injective map of ordered sets.
A typical example $R_{\mathrm{row}} = R_{\mathrm{row},\lam/\zeta} $ of such maps
 (called the {\it row reading} of $\lam/\zeta$)
is given by 
$$
R_{\mathrm{row}} (i,j) = |\{ (i^{\prime}, j^{\prime}) \in D(\lam/\zeta)\,|\, i < i^{\prime}\,\, 
 \text{or}\,\, i = i^{\prime},\,  j^{\prime} \leq j\}|.
$$
For example, if $\lam/\zeta = (5,5,4,2,1)/(3,2,1)$, then $R_{\mathrm{row}} $ is expressed as

\begin{equation}
\label{ExampleRrow}
\setlength{\unitlength}{.8mm}
\begin{picture}(5,35)
 \put(5,20){\makebox(4.9,4.9){$R_{\mathrm{row}} \,\,\,=$}}
\end{picture}
\qquad\quad
\setlength{\unitlength}{0.8mm}
\begin{picture}(25,35)
 \put(0,5){\line(0,1){25}}
 \put(0,30){\line(1,0){25}}
 \put(0,5){\framebox(4.9,4.9){$1$}}
 \put(0,10){\framebox(4.9,4.9){$2$}}
 \put(5,10){\framebox(4.9,4.9){$3$}}
 \put(5,15){\framebox(4.9,4.9){$4$}}
 \put(10,15){\framebox(4.9,4.9){$5$}}
 \put(15,15){\framebox(4.9,4.9){$6$}}
 \put(10,20){\framebox(4.9,4.9){$7$}}
 \put(15,20){\framebox(4.9,4.9){$8$}}
 \put(20,20){\framebox(4.9,4.9){$9$}}
 \put(15,25){\framebox(4.9,4.9){$10$}}
 \put(20,25){\framebox(4.9,4.9){$11$}}
 \put(25,19){\makebox(4.9,4.9){.}}
\end{picture}
\end{equation}

We say that a partial tableau $T$ on $\nu/\eta$ is a
{\it Remmel-Whitney tableau} of type $R$ (\cite{RemmelWhitney}) 
if  it satisfies the following three conditions:
\begin{quote}
(1) The image of $T$ coincides with that of  $R$.\\
(2) For each  $(i,j), (i,j+1) \in D(\lam/\zeta)$, we have $T^{-1} (a) \leqsw T^{-1} (b)$, 
where $a = R(i,j)$ and $b = R(i,j+1)$.\\
(3) For each  $(i,j), (i+1,j) \in D(\lam/\zeta)$, we have $T^{-1} (a) \leqsw T^{-1} (b)$, 
where $a = R(i,j)$ and $b = R(i+1,j)$.
\end{quote}

The correspondence $\Psf \mapsto R \circ \Psf$ gives a bijection
from $\Pic (\nu/\eta, \lam/\zeta)$ onto the set $\mathrm{RW} (\nu/\eta;R)$ of Remmel-Whitney tableau on  $\nu/\eta$ of type $R$
(cf. \cite{FominGreene}).    

\begin{example}
\label{ExampleLL}
{\rm
Let $\lam, \mu$ and $\zeta$ be $(5,5,4,2,1)$, $(4,4,4,3,2)$ and $(3,3,2,1)$, 
respectively.
Let $R\!: D(\lam/\zeta) \to \Z_{>0}$  be  the restriction of 
\eqref{ExampleRrow} and let $T$ be as in \eqref{ExampleTEq}.
By checking six order relations including 
$$
T^{-1} (R(1,4)) = T^{-1} (10) 
= (4,3) \leqsw (3,4) = 
T^{-1} (11) = T^{-1} (R(1,5)),
$$
we see that $T$ is an element of $\mathrm{RW} ((\mu/\zeta)^{\tp} ;R)$. 
Hence, there exists a unique  $\LL \in \Pic ((\mu/\zeta)^{\tp}, \lam/\zeta)$ which satisfies
$T = R \circ \LL$. 
Explicitly, $\LL$ is given by
$\LL (a) = A, \LL (b) = B,\ldots, \LL (h) = H$, 
where $a = (4,1), b,\ldots, H \in \Z^2_{> 0}$
are as in \eqref{ExampleLLEq} below, 
since $\LL (a) = R^{-1}(T(a)) = R^{-1}(1) = A$, for example.

\begin{equation}
\label{ExampleLLEq}
\qquad
\setlength{\unitlength}{0.8mm}
\begin{picture}(15,35)
 \put(0,10){\line(0,1){20}}
 \put(0,30){\line(1,0){25}}
 \put(0,10){\framebox(4.9,4.9){$a$}}
 \put(5,10){\framebox(4.9,4.9){$c$}}
 \put(10,10){\framebox(4.9,4.9){$g$}}
 \put(10,15){\framebox(4.9,4.9){$b$}}
 \put(15,15){\framebox(4.9,4.9){$h$}}
 \put(15,20){\framebox(4.9,4.9){$e$}}
 \put(20,20){\framebox(4.9,4.9){$f$}}
 \put(20,25){\framebox(4.9,4.9){$d$}}
\end{picture}
\qquad\qquad
\setlength{\unitlength}{.8mm}
\begin{picture}(5,35)
 \put(5,20){\makebox(4.9,4.9){\scalebox{1.5}{$\overset{\LL}{\longrightarrow}$}}}
\end{picture}
\qquad\qquad
\setlength{\unitlength}{0.8mm}
\begin{picture}(25,35)
 \put(0,5){\line(0,1){25}}
 \put(0,30){\line(1,0){25}}
 \put(0,5){\framebox(4.9,4.9){$A$}}
 \put(5,10){\framebox(4.9,4.9){$B$}}
 \put(10,15){\framebox(4.9,4.9){$C$}}
 \put(15,15){\framebox(4.9,4.9){$D$}}
 \put(15,20){\framebox(4.9,4.9){$E$}}
 \put(20,20){\framebox(4.9,4.9){$F$}}
 \put(15,25){\framebox(4.9,4.9){$G$}}
 \put(20,25){\framebox(4.9,4.9){$H$}}
 \put(25,15){\makebox(4.9,4.9){.}}
\end{picture}
\end{equation}

}
\end{example}

Let  $\Psf \in \Pic (\nu/\eta, \lam/\zeta)$ be a picture and let $\zz$ be an element
of $\bICC (\lam/\zeta)$.
Let $R\!:(D(\lam/\zeta)\amalg \{ \zz \},\leqsw)  \to \Z_{>0}$ be an injective map of ordered sets. 
Then the bumping route 
$\xx_l, \xx_{l-1}, \ldots, \xx_r$ of $(R \circ \Psf,R(\zz))$
does not depend on the choice of $R$.
In fact, the route is characterized by the following three conditions(cf.\cite{Leeuwen}, page 338):

\begin{quote}
{\rm (1)} $\Psf(\xx_i) \,\, (i  > r)$ is the predecessor of $\yy_{i}$ in the totally ordered set $(\Psf (D_i (\nu/\eta)) \amalg \{\yy_i \},\leqsw)$,
where  $\yy_j = \Psf (\xx_{j+1})$ for  $r \leq j<l$ and $\yy_l = \zz$.\\
{\rm (2)}  If $r > 0$, then $\yy_r$ is the minimal element of  
$(\Psf (D_r (\nu/\eta)) \amalg \{\yy_r \},\leqsw)$.\\
{\rm (3)} Either $r > 0$ and $\xx_r = (r, \eta_{r})$,  or $r = 0$
and $\xx_0 = (0,\nu_1)$. 
\end{quote}

We call $\xx_l, \xx_{l-1}, \ldots, \xx_r$ 
and $\uu := \xx_r \in \bICC (\nu/\eta)$ the {\it bumping route} and 
the {\it bumping destination} of $(\Psf,\zz)$, respectively.
Also, we write
$\displaystyle{\uu\, \frac{\,\,\,\Psf\,\,\,}{}\,\zz}$.
When $\zz \in \mathrm{ICC} (\lam/\zeta)$ and $\uu  \in \mathrm{ICC} (\nu/\eta)$, we say that 
$\zz$ is 
an {\it addable cocorner} of $\Psf$.
In this case, we define the {\it row insertion}  
$E_{\zz} \Psf$ to be the unique picture from $\nu/(\eta \setminus  \uu)$ onto 
$\lam/(\zeta \setminus  \zz) $ 
satisfying the condition
\begin{equation}
\label{ERPi}
 E_{R(\zz)} (R \circ \Psf) = R \circ (E_{\zz} \Psf)
\end{equation}
 (cf. \cite{Zelevinsky}, \cite{Leeuwen}).
More explicitly, the picture $E_{\zz}  \Psf$ is given by  
$(E_{\zz} \Psf) (\xx_l) = \zz $, 
$(E_{\zz} \Psf) (\xx_i) = \Psf(\xx_{i + 1}) \,\,\,\, $ $(l>  i \geq r)$ 
and 
$(E_{\zz} \Psf)(\yy) = \Psf(\yy) \,\,\,\,  (\yy \ne \xx_l, \xx_{l-1}, \ldots, \xx_r)$.\\

Let  $\Psf \in \Pic (\nu/\eta, \lam/\zeta)$ be a picture and let $\vv$ be an inner corner of  
of $\neta$.
Then there exists at most one pair $(\Psf^{\prime}, \ww)$ such that 
$\Psf = E_{\ww} \Psf^{\prime}$ and that $\vv$ is the bumping destination of $(\Psf^{\prime}, \ww)$.
When there exists such a pair, 
we say that $\vv$ is a  {\it removable corner} of $\Psf$ and we write
$\Psf^{\prime} = F_{\vv} \Psf$ and $\displaystyle{\vv\, \frac{}{\,\,\,\Psf\,\,\,}\, \ww}$. 
We note that $\ww$ is an inner corner of $\lam/\zeta$.
Let $R\!:(D(\lam/\zeta),\leqsw)  \to \Z_{>0}$ be an injective map of ordered sets. 
Then,  $\vv$ is a {\it removable corner} of $\Psf$ if  and only if  it is a removable corner of 
$R \circ \Psf$.  Moreover, we have 
\begin{gather}
\label{FRPi}
 F_{\vv} (R \circ \Psf) = R \circ (F_{\vv} \Psf),\\
\label{v-Rw}
 \vv\, \frac{}{\,\,\,R \circ \Psf\,\,\,}\, R(\ww).
\end{gather}

\begin{example}
\label{ExampleELL}
{\rm
Let $\lam, \mu, \zeta$ and $\LL$ be as in 
Example \ref{ExampleLL}.
Let $R\!: D(\lam/\zeta)\amalg \{ (2,3) \} \to \Z_{>0}$  be  the restriction of 
\eqref{ExampleRrow}. Then $R(2,3) = 7$ and $T = R \circ \LL$ is given by \eqref{ExampleTEq}.
Hence by  Example \ref{ExampleT} and \eqref{ERPi}, we have $\displaystyle{(2,3)\, \frac{\,\,\,\LL \,\,\,}{}\, (2,3)}$ and

\begin{equation}
\label{ExampleELLEq}
\qquad
\setlength{\unitlength}{0.8mm}
\begin{picture}(15,35)
 \put(0,10){\line(0,1){20}}
 \put(0,30){\line(1,0){25}}
 \put(0,10){\framebox(4.9,4.9){$a$}}
 \put(5,10){\framebox(4.9,4.9){$e$}}
 \put(10,10){\framebox(4.9,4.9){$h$}}
 \put(10,15){\framebox(4.9,4.9){$c$}}
 \put(15,15){\framebox(4.9,4.9){$i$}}
 \put(10,20){\framebox(4.9,4.9){$b$}}
 \put(15,20){\framebox(4.9,4.9){$f$}}
 \put(20,20){\framebox(4.9,4.9){$g$}}
 \put(20,25){\framebox(4.9,4.9){$d$}}
\end{picture}
\qquad\qquad
\setlength{\unitlength}{.8mm}
\begin{picture}(5,35)
 \put(5,20){\makebox(4.9,4.9){\scalebox{1.5}{$\overset{E_{(2,3)} \LL}{\longrightarrow}$}}}
\end{picture}
\qquad\qquad
\setlength{\unitlength}{0.8mm}
\begin{picture}(25,35)
 \put(0,5){\line(0,1){25}}
 \put(0,30){\line(1,0){25}}
 \put(0,5){\framebox(4.9,4.9){$A$}}
 \put(5,10){\framebox(4.9,4.9){$B$}}
 \put(10,15){\framebox(4.9,4.9){$C$}}
 \put(15,15){\framebox(4.9,4.9){$D$}}
 \put(10,20){\framebox(4.9,4.9){$E$}}
 \put(15,20){\framebox(4.9,4.9){$F$}}
 \put(20,20){\framebox(4.9,4.9){$G$}}
 \put(15,25){\framebox(4.9,4.9){$H$}}
 \put(20,25){\framebox(4.9,4.9){$I$}}
 \put(25,19){\makebox(4.9,4.9){.}}
\end{picture}
\end{equation}
Here the notation is the same as in \eqref{ExampleLLEq}, that is, 
$(E_{(2,3)} \LL)(a) = A$, for example.
Similarly, by the last assertions of Example \ref{ExampleT} and \eqref{FRPi}, we have
$\displaystyle{(3,3)\, \frac{}{\,\,\,\LL\,\,\,}\, (3,3)}$ and 

\begin{equation}
\label{ExampleFLLEq}
\qquad
\setlength{\unitlength}{0.8mm}
\begin{picture}(15,35)
 \put(0,10){\line(0,1){20}}
 \put(0,30){\line(1,0){25}}
 \put(0,10){\framebox(4.9,4.9){$a$}}
 \put(5,10){\framebox(4.9,4.9){$b$}}
 \put(10,10){\framebox(4.9,4.9){$f$}}
 \put(15,15){\framebox(4.9,4.9){$g$}}
 \put(15,20){\framebox(4.9,4.9){$d$}}
 \put(20,20){\framebox(4.9,4.9){$e$}}
 \put(20,25){\framebox(4.9,4.9){$c$}}
\end{picture}
\qquad\qquad
\setlength{\unitlength}{.8mm}
\begin{picture}(5,35)
 \put(5,20){\makebox(4.9,4.9){\scalebox{1.5}{$\overset{F_{(3,3)} \LL}{\longrightarrow}$}}}
\end{picture}
\qquad\qquad
\setlength{\unitlength}{0.8mm}
\begin{picture}(25,35)
 \put(0,5){\line(0,1){25}}
 \put(0,30){\line(1,0){25}}
 \put(0,5){\framebox(4.9,4.9){$A$}}
 \put(5,10){\framebox(4.9,4.9){$B$}}
 \put(15,15){\framebox(4.9,4.9){$C$}}
 \put(15,20){\framebox(4.9,4.9){$D$}}
 \put(20,20){\framebox(4.9,4.9){$E$}}
 \put(15,25){\framebox(4.9,4.9){$F$}}
 \put(20,25){\framebox(4.9,4.9){$G$}}
 \put(25,15){\makebox(4.9,4.9){.}}
 \put(5,15){\line(1,0){15}}
\end{picture}
\end{equation}

}
\end{example}

\begin{lemma}
\label{F=E^-1}
 $\mathrm{(1)}$ 
Let $\zz$ be an addable cocorner of $\Psf$ and set 
 $\displaystyle{\uu\, \frac{\,\,\,\Psf\,\,\,}{}\, \zz}$. 
Then $\uu$ is a removable corner of $E_{\zz} \Psf$, 
which satisfies $\displaystyle{\uu\, \frac{}{\,\,\, E_{\zz} \Psf\,\,\,}\, \zz}$. 
Moreover, we have 
\begin{equation}
\label{FEPi}
F_{\uu} (E_{\zz}\Psf )= \Psf.
\end{equation}
$\mathrm{(2)}$ 
Let $\vv$ be a removable corner of $\Psf$ and set 
 $\displaystyle{\vv\, \frac{}{\,\,\,\Psf\,\,\,}\, \ww}$. 
Then $\ww$ is an addable cocorner of $F_{\vv} \Psf$, 
which satisfies $\displaystyle{\vv\, \frac{\,\,\, F_{\vv} \Psf\,\,\,}{}\, \ww}$.
Moreover, we have 
\begin{equation}
\label{EFPi}
E_{\ww} (F_{\vv}\Psf )= \Psf.
\end{equation}
\end{lemma}

As an immediate consequence of Lemma \ref{BumpLemTab}, 
we have the following:
\begin{lemma}
\label{BumpLemPict}
For each $\Psf \in \Pic (\nu/\eta, \lam/\zeta) $, we have the following:\\
 $\mathrm{(1)}$ 
Let $\zz$ and $\zzp$ be elements of $\bICC (\lam/\zeta)$ such that
$\zz \leq_{\swarrow} \zzp$. Define $\uu, \uup \in \bICC(\mu/\eta)$ by  
$\displaystyle{\uu\, \frac{\,\,\,\Psf\,\,\,}{}\, \zz}$, $\displaystyle{\uup\, \frac{\,\,\,\Psf\,\,\,}{}\, \zzp}$. 
Then we have $\uu \leq_{\swarrow} \uup$.
 \\
%$\mathrm{(2)}$ Let $\vv$ and $\vvp$ be elements of $\IC (\mu/\eta)$
%such that 
%$\vv \leq_{\swarrow} \vvp$ and that 
%$\vvp$ is a removable corner of $\Psf$.
%Then $\vv$ is also a removable corner of $\Psf$.
%Moreover,  we have $\ww \leq_{\swarrow} \wwp$, where
%$\displaystyle{\vv\, \frac{}{\,\,\,\Psf\,\,\,}\, \ww}$, $\displaystyle{\vvp\, \frac{}{\,\,\,\Psf\,\,\,}\, \wwp}$. 
%
%
$\mathrm{(2)}$ Let $\vv$ and $\vvp$ be removable corners of $\Psf$ such that
$\vv \leq_{\swarrow} \vvp$.
Define $\ww, \wwp \in \IC (\lz)$ by
$\displaystyle{\vv\, \frac{}{\,\,\,\Psf\,\,\,}\, \ww}$, $\displaystyle{\vvp\, \frac{}{\,\,\,\Psf\,\,\,}\, \wwp}$. 
Then, we have $\ww \leq_{\swarrow} \wwp$.\\
$\mathrm{(3)}$ 
 Let $\zz$ be an addable cocorner of $\Psf$ and let $\zzp$ be 
an inner cocorner of $\lam/(\zeta \setminus \zz)$.
 Define points $\uu$ and  $\uup$ by
$\displaystyle{\uu\, \frac{\,\,\,\Psf\,\,\,}\, \zz}{}$ and $\displaystyle{\uup\, \frac{\,\,\,E_{\zz} \Psf\,\,\,}{}\, \zzp}$,
respectively. 
Then we have
$\zz \lsw \zzp \Rightarrow \uu \lsw \uup$ and 
$\zzp \lsw \zz \Rightarrow \uup \lsw \uu$.\\
%
%
%$\mathrm{(4)}$
%Let $\vv$ be a removable corner of $\Psf$ with
%$\displaystyle{\vv\, \frac{}{\,\,\,\Psf\,\,\,}\, \ww}$
%and let  $\vvp$ be an inner corner of 
%$\nu/(\eta \amalg \vv)$.
%If $\vvp \lsw \vv$, then $\vvp$ is also removable and satisfies $\wwp \lsw \ww$, where 
%$\displaystyle{\vvp\, \frac{}{\,\,\,F_{\vv} \Psf\,\,\,}\, \wwp}$.
%If $\vv \lsw \vvp$  and $\vvp$ is also removable, then  $\ww \lsw \wwp$, where 
%$\displaystyle{\vvp\, \frac{}{\,\,\,F_{\vv} \Psf\,\,\,}\, \wwp}$.
%
%
$\mathrm{(4)}$
Let $\vv$ be a removable corner of $\Psf$ and 
let $\vvp$ be a removable corner of $F_{\vv} \Psf$.
Define points $\ww$ and $\wwp$ by 
$\displaystyle{\vv\, \frac{}{\,\,\,\Psf\,\,\,}\, \ww}$ and 
$\displaystyle{\vvp\, \frac{}{\,\,\,F_{\vv} \Psf\,\,\,}\, \wwp}$, respectively. 
Then, we have 
$\vv \lsw \vvp \Rightarrow \ww \lsw \wwp$ and $\vvp \lsw \vv \Rightarrow \wwp \lsw \ww$.\\
%
%
%\par\noindent 
%$\mathrm{(5)}$
%Let $\uu$ be an inner cocorner of $\nu/\eta$
%and let $\vv$  be an inner corner of $\nu/\eta$.
%Suppose that there exists an addable cocorner $\zz$ of $\Psf$ such that
%$\displaystyle{\uu\, \frac{\,\,\,\Psf\,\,\,}{}\, \zz}$.
%If $\vv \lsw \uu$, then 
% $\vv$ is a removable corner of $\Psf$ and satisfies  $\ww \lsw \zz$, where
%$\displaystyle{\vv\, \frac{}{\,\,\,\Psf\,\,\,}\, \ww}$.
%If $\uu \lsw \vv$ and 
% $\vv$ is a removable corner of $\Psf$ with $\displaystyle{\vv\, \frac{}{\,\,\,\Psf\,\,\,}\, \ww}$, 
%then $\zz \lsw \ww$.
%
%
$\mathrm{(5)}$
Let $\zz$ be an element of $\bICC (\lz)$ and
define $\uu \in \bICC (\nu/\eta)$ by 
$\displaystyle{\uu\, \frac{\,\,\,\Psf\,\,\,}{}\, \zz}$.
Let $\vv$ be a removable corner of $\Psf$ such that $\uu \lsw \vv$. 
Then, we have $\zz \lsw \ww$, where $\displaystyle{\vv\, \frac{}{\,\,\,\Psf\,\,\,}\, \ww}$.\\
$\mathrm{(6)}$
Let $\zz$ and $\uu$ be as in Part (5).
If a cocorner $\vv$ of $\eta$ satisfies $\vv \lsw \uu$, then
$\vv$ is a removable corner of $\Psf$.
Moreover, we have $\ww \lsw \zz$, where  
$\displaystyle{\vv\, \frac{}{\,\,\,\Psf\,\,\,}\, \ww}$.
\end{lemma}

\subsection{Balanced corners}
\label{Balanced corners}

Let $\lam$ and $\mu$ be partitions of $n$ and let $\zeta$ be a partition of
$n - m$,  which satisfies both $\zeta \subseteq \lam$  and $\zeta \subseteq \mu$.
We define a set $\PW (\lambda, \mu; \zeta)$ by \eqref{DefPElmz} and
call its element a {\it picture of type} $(\lambda, \mu; \zeta)$.
Also, we define a set 
$\PW_m (\lambda, \mu)$ by
\eqref{DefPEm}
and call its element a {\it picture of type} 
 $(\lambda, \mu)$ with size $m$.
For convenience, we set $\PW_0 (\lambda, \mu) = \emptyset$ if
$\lam \ne \mu$ and
$$
 \PW_0 (\lam, \lam) = \PW (\lam,\lam;\lam) = \{ \id_{\emptyset} \}.
$$
Let $\LL$ be an element of $\PW_m (\lambda, \mu)$ and let
$\vv$ be its removable corner.
We say that  $\vv$   is a 
{\it balanced corner} of $\LL$ if 
$$\displaystyle{ \vv\, \frac{}{\,\,\,\LL\,\,\,}\,\vv^{\tp}}.$$
%$$\displaystyle{ {\boldsymbol v}\, \frac{}{\,\,\,\LL\,\,\,}\,{\boldsymbol v}^{\tp}}.$$
%
Also, 
we say that an addable cocorner $\zz$ of $\VV \in \PW_m (\lambda, \mu)$
is a {\it balanced cocorner} if 
$\displaystyle{  \zz^{\tp}\, \frac{\,\,\,{ \VV}{}\,\,\,}\, \zz}.$

\begin{example}
\label{ExampleBC}
{\rm 
Let  $\LL \in \PW ((5,5,4,2,1), (4,4,4,3,2); (3,3,2,1))$ be as in \eqref{ExampleLLEq}.
By Example  \ref{ExampleELL},  $(3,3)$ is a balanced corner of $\LL$,
while $(2,3)$ is not a balanced cocorner of $\LL$.

}
\end{example}

\begin{lemma}
\label{BCC}
For each $\LL, \VV \in \PW (\lambda, \mu; \zeta)$, we have the following:\\
$\mathrm{(1)}$ 
The picture $\LL$ has at most one  
balanced corner. Also, it  has at most one balanced cocorner.\\
$\mathrm{(2)}$ 
If $\VV$ has a 
balanced cocorner $\zz$, then $E_{\zz} \VV$ has no balanced cocorners.
If $\LL$ has a 
balanced corner $\vv$, then $F_{\vv} \LL$ has no balanced corners. \\
$\mathrm{(3)}$ 
If $\VV$ has a 
balanced cocorner,  then $\VV$ has no balanced corners.
If $\LL$ has a 
balanced corner,  then $\LL$ has no balanced cocorners.
\end{lemma}

\noindent
{\it Proof.} 
Let $\zz$ and $\zz^{\prime}$ be balanced cocorners of $\LL$.
Since $(\mathrm{ICC} (\lam/\zeta), \leqsw)$ is  totally ordered, 
we may assume $\zz \leqsw \zz^{\prime}$.
Then, by Lemma \ref{BumpLemPict} (1),  we have $\zz^{\tp} \leqsw (\zz^{\prime})^{\tp}$.
Hence, we have $ \zz^{\prime} \leqsw \zz$ by \eqref{Tenchi&Leqsw}.  
This proves the second statement of 
Part (1).
Similarly, the first statement of Part  (1) follows from 
Lemma \ref{BumpLemPict} (2), Part (2) follows from Lemma \ref{BumpLemPict} (3), (4)
and Part (3) follows from Lemma \ref{BumpLemPict} (5),(6). $\hfill \square$

\begin{lemma}
\label{DlEmpty}
Suppose that $(\mu/\zeta)^{\tp}$ satisfies
$D_l ((\mu/\zeta)^{\tp}) = \emptyset$, 
where $l = l((\mu/\zeta)^{\tp})$.
Then, for each picture $\VV \in  \PW (\lam, \mu; \zeta)$, 
$((\zeta^{\tp})_l, l)$ %$=$  $\max (\mathrm{ICC} (\lam/\zeta),$  $\leqsw)$
is a balanced cocorner of $\VV$.
\end{lemma}

\noindent
\noindent
{\it Proof.}\,\, Since $(\zeta^{\tp})_l = (\mu^{\tp})_l > 0$, 
$\uu := (l, (\zeta^{\tp})_l)$ is a corner of $\zeta^{\tp}$. 
On the other hand, by the definition of the backward row insertion, we have 
$\displaystyle{ \uu\, \frac{\,\,\,\VV\,\,\,}{}\,\zz}$   
for every $\zz \in \mathrm{ICC} (\lam/\zeta)$.
Hence $\uu^{\tp}$ is a balanced cocorner of $\VV$.  $\hfill\square$

\begin{lemma}
\label{exactness}
If $\LL$ does not have a 
balanced cocorner,  then $\LL$ has a balanced corner.
\end{lemma}

We will give a proof of the lemma above in Sect.
\ref{The exactness}.

\subsection{Picture of hook shape}
\label{Picture of the hook}

Let $\VV$ be an element of $\PW(\lam, \mu;\zeta)$. 
We say that $\VV$ is a 
 {\it picture of hook shape} of type 
$(\lam, \mu;\zeta)$ if it has a (necessarily unique) balanced cocorner.
Let $\PH(\lam, \mu;\zeta)$ be the set of pictures of hook shape of type 
$(\lam, \mu;\zeta)$ and let $\PH_m(\lam, \mu)$ be the set
$\coprod_{\zeta \in \PP_{n-m}} \PH (\lam, \mu;\zeta)$.
By Lemma \ref{BCC} (3) and  Lemma \ref{exactness}, we have the following decomposition:
\begin{equation}
\label{DecompPW}
\PW_m (\lam, \mu) = 
\PH_m (\lam, \mu) 
\amalg
\PH^c_m (\lam, \mu),  
\end{equation}
where $\PH^c_m (\lam, \mu)$ denotes the set of pictures of type
$(\lam, \mu)$ with size $m$, which have balanced corners. 
Since $\PH^c_0 (\lam, \mu) = \PH_n (\lam, \mu) = \emptyset$, we have in particular
\begin{equation}
\label{PW0=PF0}
\PW_0 (\lam, \mu) = \PH_0 (\lam, \mu), 
\quad 
\PW_n (\lam, \mu) = \PH_n^c (\lam, \mu). 
\end{equation}
Given  $\VV\in \PH_m (\lam,\mu)$ and $\LL \in \PH^c_m (\lam,\mu)$,   we define 
pictures $E \VV$ %\in PW_{m+1}(\lam,\mu)$ 
and $F \LL$ % \in  PW_{m-1}(\lam,\mu)$ 
by 
$E \VV = E_{\zz} \VV$ and $F \LL = F_{\vv} \LL$ respectively, where $\zz$ denotes 
the unique balanced cocorner of $\VV$ and $\vv$ denotes 
the unique balanced corner of $\LL$.
By Lemma \ref{BCC} (2) and \eqref{DecompPW}, we have
$E \VV \in \PH^c_{m+1} (\lam,\mu)$
and
$F \LL \in \PH_{m-1} (\lam,\mu)$.
Hence,  by Lemma \ref{F=E^-1}, the operator $E$ gives a bijection 
\begin{equation}
\label{PHm=PHcm+1}
\PH_m (\lam, \mu) \cong  
\PH_{m + 1}^c (\lam, \mu) 
\end{equation}
whose inverse is given by $F$.

%\begin{defin}
% $E \LL = E_{\zz} \LL, \LL = F_{\uu}  \LL$
%If $\LL$ does not have a balanced cocorner,  then we set 
%$E \LL = 0 \in \C \PW_{m+1} (\lam,\mu)$.   
%\end{defin}

%For each $\LL \in \PW (\lambda, \mu; \zeta)$,  we have
%$E (E \LL) = 0$ and $F (F \LL) = 0$.

\begin{theorem}
For each $\lam \in \Pn$ and $0 \leq m < n$, we have

\begin{equation}
 [L(\lam) \otimes L(n-m, 1^m)]
 =
 \sum_{\mu \in \Pn} |\PH_m (\lam, \mu)| \cdot [L(\mu)] .
\end{equation}
\end{theorem}
\medskip\medskip
\par\noindent
{\it Proof}.
By \eqref{DecompPW},  \eqref{PW0=PF0} and \eqref{PHm=PHcm+1}, 
the integers $h_m^{\prime} :=  |\PH_m (\lam, \mu)|$ satisfy  
$w_0 = h^{\prime}_0$, $w_m = h^{\prime}_{m-1} + h^{\prime}_m\,\, (0< m <n)$ and $w_n = h^{\prime}_{n-1}$. 
Solving these equations, we get $h^{\prime}_m = \sum_{i = 0}^m (-1)^{m - i} w_i$.
Therefore, the theorem follows from Part (2) of Theorem \ref{L(lam)Lam}.
$\hfill \square$

\begin{example}
\label{ExamplePH}
{\rm
For $\lam = (5,3,1,1)$ and $\mu = (4,3,3)$,
the set $\PW_6 (\lam,\mu)$ consists of 7 elements
$\LL_i$ ($i = 1,\ldots,7$), where the corresponding Remmel-Whitney tableau
$T_i = R_{row}\circ \LL_i$ are given by

\begin{multline}
%\label{}
%
\setlength{\unitlength}{.8mm}
\begin{picture}(5,35)
 \put(5,20){\makebox(4.9,4.9){$T_1 \,\,\,=$}}
\end{picture}
\qquad
\setlength{\unitlength}{0.8mm}
\begin{picture}(25,35)
 \put(0,5){\line(0,1){20}}
 \put(0,25){\line(1,0){15}}
 \put(0,5){\framebox(4.9,4.9){$3$}}
 \put(5,10){\framebox(4.9,4.9){$5$}}
 \put(10,10){\framebox(4.9,4.9){$6$}}
 \put(5,15){\framebox(4.9,4.9){$2$}}
 \put(10,15){\framebox(4.9,4.9){$4$}}
 \put(10,20){\framebox(4.9,4.9){$1$}}
 \put(15,5){\makebox(4.9,4.9){,}}
\end{picture}
\setlength{\unitlength}{.8mm}
\begin{picture}(5,35)
 \put(5,20){\makebox(4.9,4.9){$T_2 \,\,\,=$}}
\end{picture}
\qquad
\setlength{\unitlength}{0.8mm}
\begin{picture}(25,35)
 \put(0,5){\line(0,1){20}}
 \put(0,25){\line(1,0){15}}
 \put(0,5){\framebox(4.9,4.9){$2$}}
 \put(5,10){\framebox(4.9,4.9){$5$}}
 \put(10,10){\framebox(4.9,4.9){$6$}}
 \put(5,15){\framebox(4.9,4.9){$3$}}
 \put(10,15){\framebox(4.9,4.9){$4$}}
 \put(10,20){\framebox(4.9,4.9){$1$}}
 \put(15,5){\makebox(4.9,4.9){,}}
\end{picture}
\setlength{\unitlength}{.8mm}
\begin{picture}(5,35)
 \put(5,20){\makebox(4.9,4.9){$T_3 \,\,\,=$}}
\end{picture}
\qquad
\setlength{\unitlength}{0.8mm}
\begin{picture}(25,35)
 \put(0,5){\line(0,1){20}}
 \put(0,25){\line(1,0){15}}
 \put(0,5){\framebox(4.9,4.9){$5$}}
 \put(5,10){\framebox(4.9,4.9){$3$}}
 \put(10,10){\framebox(4.9,4.9){$6$}}
 \put(5,15){\framebox(4.9,4.9){$2$}}
 \put(10,15){\framebox(4.9,4.9){$4$}}
 \put(10,20){\framebox(4.9,4.9){$1$}}
 \put(15,5){\makebox(4.9,4.9){,}}
\end{picture}
\setlength{\unitlength}{.8mm}
\begin{picture}(5,35)
 \put(5,20){\makebox(4.9,4.9){$T_4 \,\,\,=$}}
\end{picture}
\qquad\quad
\setlength{\unitlength}{0.8mm}
\begin{picture}(25,35)
 \put(0,5){\line(0,1){20}}
 \put(0,25){\line(1,0){15}}
 \put(0,5){\framebox(4.9,4.9){$4$}}
 \put(0,10){\framebox(4.9,4.9){$3$}}
 \put(5,10){\framebox(4.9,4.9){$5$}}
 \put(10,10){\framebox(4.9,4.9){$6$}}
 \put(10,15){\framebox(4.9,4.9){$2$}}
 \put(10,20){\framebox(4.9,4.9){$1$}}
 \put(15,5){\makebox(4.9,4.9){,}}
\end{picture}\\
\setlength{\unitlength}{.8mm}
\begin{picture}(5,35)
 \put(5,20){\makebox(4.9,4.9){$T_5 \,\,\,=$}}
\end{picture}
\qquad\quad
\setlength{\unitlength}{0.8mm}
\begin{picture}(25,35)
 \put(0,5){\line(0,1){20}}
 \put(0,25){\line(1,0){15}}
 \put(0,5){\framebox(4.9,4.9){$4$}}
 \put(0,10){\framebox(4.9,4.9){$2$}}
 \put(5,10){\framebox(4.9,4.9){$5$}}
 \put(10,10){\framebox(4.9,4.9){$6$}}
 \put(10,15){\framebox(4.9,4.9){$3$}}
 \put(10,20){\framebox(4.9,4.9){$1$}}
 \put(15,5){\makebox(4.9,4.9){,}}
\end{picture}
\setlength{\unitlength}{.8mm}
\begin{picture}(5,35)
 \put(5,20){\makebox(4.9,4.9){$T_6 \,\,\,=$}}
\end{picture}
\qquad\quad
\setlength{\unitlength}{0.8mm}
\begin{picture}(25,35)
 \put(0,5){\line(0,1){20}}
 \put(15,20){\line(0,1){5}}
 \put(0,25){\line(1,0){15}}
 \put(0,5){\framebox(4.9,4.9){$4$}}
 \put(0,10){\framebox(4.9,4.9){$2$}}
 \put(5,10){\framebox(4.9,4.9){$5$}}
 \put(10,10){\framebox(4.9,4.9){$6$}}
 \put(5,15){\framebox(4.9,4.9){$1$}}
 \put(10,15){\framebox(4.9,4.9){$3$}}
 \put(15,5){\makebox(4.9,4.9){,}}
\end{picture}
\setlength{\unitlength}{.8mm}
\begin{picture}(5,35)
 \put(5,20){\makebox(4.9,4.9){$T_7 \,\,\,=$}}
\end{picture}
\qquad\quad
\setlength{\unitlength}{0.8mm}
\begin{picture}(25,35)
 \put(0,5){\line(0,1){20}}
 \put(15,20){\line(0,1){5}}
 \put(0,25){\line(1,0){15}}
 \put(0,5){\framebox(4.9,4.9){$4$}}
 \put(0,10){\framebox(4.9,4.9){$1$}}
 \put(5,10){\framebox(4.9,4.9){$5$}}
 \put(10,10){\framebox(4.9,4.9){$6$}}
 \put(5,15){\framebox(4.9,4.9){$2$}}
 \put(10,15){\framebox(4.9,4.9){$3$}}
 \put(15,5){\makebox(4.9,4.9){.}}
\end{picture}
\end{multline}
The point $(1,3)$ is a balanced cocorner of both $\LL_1$ and $\LL_2$.
While the points $(1,4)$ and $(1,3)$ are balanced corners of $\LL_3$ and $\LL_i$ $(i \geq 4)$, respectively.
Hence, the multiplicity of $L(4,3,3)$ in $L(5,3,1,1)\otimes L(4,1^6)$ is $2$
and that of $L(4,3,3)$ in  $L(5,3,1,1)\otimes L(5,1^5)$ is $5$.
}

\end{example}

\begin{example}
\label{Hook^3}
{\rm (} {\rm \cite{Rosas}, Theorem 3} {\rm )}
Let $\lambda$ be $(n-e,1^e)$ and let $\mu$ be $(n-f,1^f)$, 
where $2e, 2f \leq n$ and $g: = f - e \geq 0$.
Then $\PH_{g+i} (\lambda,\mu)\ne\emptyset$
if and only if either $e + f < n$ and $0 \leq i \leq 2e$,
or $e + f = n$ and $0 \leq i \leq n-2$.
In this case, we have
$\PH_{g+i} (\lambda,\mu) = \{ \VV_i\}$, where
$T_i = R_{row}\circ \VV_i$ are given by  
\begin{equation}
\label{PicHook^3}
\setlength{\unitlength}{0.8mm}
\begin{picture}(15,10)
 \put(07,40){\makebox(4.9,4.9){$T_{2k-1} =$}}
\end{picture}\qquad
\setlength{\unitlength}{0.8mm}
\begin{picture}(25,55)
 \put(33,50){\mbox{$f$}}
 \qbezier(5,46)(33,50)(62,46)
 \put(-6.8,25){\mbox{$v$}}
 \qbezier(-1.5,5)(-5.5,25)(-1.5,45)
 \put(0,5){\line(0,1){40}}
 \put(5,5){\line(0,1){35}}
 \put(5,40){\line(1,0){40}}
 \put(0,45){\line(1,0){45}}
 \put(0,20){\framebox(4.9,4.9){$2$}}
 \put(0,13){\mbox{$\,\,\vdots$}}
 \put(0,5){\framebox(4.9,4.9){$k$}}
 \put(15,40){\framebox(4.9,4.9){$1$}}
 \put(20,40){\framebox(12.8,4.9){$k+1$}}
 \put(34,41){\mbox{$\cdots$}}
 \put(42,40){\framebox(19.8,4.9){$g\!+\!2k\!-\!1$}}
\end{picture}
\qquad\qquad\qquad\qquad\qquad
\setlength{\unitlength}{0.8mm}
\begin{picture}(15,10)
 \put(0,40){\makebox(4.9,4.9){$\,\,\,,\quad T_{2k} =$\qquad\,\,\,}}
\end{picture}
\setlength{\unitlength}{0.8mm}
\begin{picture}(15,55)
 \put(28,50){\mbox{$f$}}
 \qbezier(5,46)(28,50)(52,46)
 \put(-6.8,25){\mbox{$v$}}
 \qbezier(-1.5,5)(-5.5,25)(-1.5,45)
 \put(0,5){\line(0,1){40}}
 \put(5,5){\line(0,1){35}}
 \put(5,40){\line(1,0){35}}
 \put(0,45){\line(1,0){40}}
 \put(0,20){\framebox(4.9,4.9){$1$}}
 \put(0,13){\mbox{$\,\,\vdots$}}
 \put(0,5){\framebox(4.9,4.9){$k$}}
 \put(15,40){\framebox(12.8,4.9){$k+1$}}
 \put(29,41){\mbox{$\cdots$}}
 \put(37,40){\framebox(14.8,4.9){$g+2k$}}
 \put(47,39){\makebox(14.8,4.9){$.$}}
\end{picture}
\setlength{\unitlength}{0.8mm}
\begin{picture}(40,20)
 \put(10,15){\makebox(4.9,4.9){$\,$}}
\end{picture}
%\setlength{\unitlength}{0.8mm}
%\begin{picture}(40,20)
% \put(-1,39){\makebox(4.9,4.9){$.$}}
%\end{picture}
\end{equation}

\end{example}

\subsection{The exactness}
\label{The exactness}

In this section, we will give a proof of Lemma \ref{exactness}.
Throughout this subsection, we fix a picture $\LL$ of type 
 $(\lambda, \mu; \zeta)$, which has no balanced cocorners.
For each $\zz \in \bICC (\lam/\zeta)$, define $\uu(\zz) \in \bICC((\mu/\zeta)^{\tp})$ by
$\displaystyle{ \uu(\zz)\, \frac{\,\,\,\LL\,\,\,}{}\,\zz}$.
We also use the following notations:
\begin{gather*}
 l = l ((\mu/\zeta)^{\tp}), \\
\quad
\uu_{--}:= {\min}_{\leqsw} \bICC ((\mu/\zeta)^{\tp}),
\quad
\uu_{++}:= {\max}_{\leqsw} \bICC ((\mu/\zeta)^{\tp}),
\\
\zz_{--} := {\min}_{\leqsw} \bICC (\lam/\zeta),
\quad
\zz_{++} := {\max}_{\leqsw} \bICC (\lam/\zeta),
\\
C_{-} := \{  \zz \in \bICC (\lam/\zeta))\,|\,  \zz \lsw \uu (\zz)^{\tp} \},\\
%\quad
C_{+} := \{  \zz \in \bICC (\lam/\zeta))\,|\,  \uu (\zz)^{\tp} \lsw \zz \}.
\end{gather*}

\begin{example}
{\rm Let $\LL$ be as in Example \ref{ExampleBC} and let $\zz_i$  
be as in the first figure of \eqref{ExampleIC}.
Then we have 
$\zz_0 \lsw  (0,4) = \uu(\zz_0)^{\tp}$, 
$\zz_1 \lsw \zz_2 \lsw \zz_3 = \uu(\zz_1)^{\tp} = \uu(\zz_2)^{\tp}$
and
$\uu(\zz_3)^{\tp} = \uu(\zz_4)^{\tp} = \zz_2 \lsw \zz_3 \lsw \zz_4$.
Hence $C_{-} = \{\zz_0, \zz_1, \zz_2 \}$ and 
$C_{+} = \{\zz_3, \zz_4\}$.

}
\end{example}

\begin{lemma}
\label{u--<uz++}
We have 
$ \uu (\zz_{--}) = \uu_{--} \lsw \uu (\zz_{++}) $.
%Moreover, if  a cocorner $\vv = (i,j)$  of $\zeta^{\tp}$ satisfies 
%$\uu_{--}  \lsw 
%$\vv \lsw \uu (\zz_{++}) $,
%then $\vv$ is a removable corner of $\LL$.   
\end{lemma}

\noindent
{\it Proof.} 
Let $R\!: (D(\lz)\amalg \{\zz_{--},\zz_{++}\},\leqsw) \to \Z_{>0}$ be an injective map of ordered sets.
Then, we have $R(\zz_{--}) < T(\xx) <R(\zz_{++})$ for every $\xx \in D((\mu/\zeta)^{\tp})$, 
where $T := R \circ \LL$.
On the other hand, by Lemma \ref{DlEmpty}, the last row  
$D_l ((\mu/\zeta)^{\tp})$ of $(\mu/\zeta)^{\tp}$ is not empty.
Hence the assertion follows immediately 
from the definition of the backward row insertion.
 $\hfill\square$

\begin{lemma}
 We have $\zz_{--} \in C_{-}$ and $\zz_{++} \in C_+$.
 In particular, $C_-, C_+ \ne \emptyset$.
\end{lemma}

\noindent
{\it Proof.} 
Suppose that $\uu (\zz_{++})$ is not extreme. Then we have 
$\uu (\zz_{++})^{\tp} \in \ICC (\lz)$.
Hence $\uu (\zz_{++})^{\tp} \leqsw \zz_{++}$.
Since $\zz_{++}$ is not a balanced cocorner even if $\zz_{++} \in \ICC (\lz)$, 
this proves $\zz_{++} \in C_+$. 
Next, suppose that $\uu (\zz_{++})$ is extreme. By Lemma \ref{u--<uz++}, we have 
$\uu (\zz_{++}) = \uu_{++} = (0,(\mu^{\tp})_1)$. Hence $\zz_{++} \in C_+$. 
The proof of $\zz_{--} \in C_-$ is similar and more easy, 
since $\uu (\zz_{--}) = \uu_{--}$.
 $\hfill\square$

In view of the lemma above, we may define
$\zz_{-} \in C_-$ and $\zz_{-} \in C_+$ by
$$
\zz_{-} := {\max}_{\leqsw} C_{-},
\quad
\zz_{+} := {\min}_{\leqsw} C_{+}.
$$

\begin{lemma}
$\mathrm{(1)}$ 
If $\zz \in \bICC (\lz)$ satisfies $\zz \leqsw \zz_{-}$, then $\zz \in C_{-}$.
\\
$\mathrm{(2)}$ 
 For each $\zz \in C_{-}$ and $\zzp \in C_+$, we have $\zz \lsw \zzp$. 
In particular, $\zz_+$ is the successor of $\zz_{-}$ in
$(\bICC (\lz), \leqsw)$.
\end{lemma}
\noindent
{\it Proof.} 
If $z \in \bICC (\lz)$ satisfies $\zz \leqsw \zz_{-}$, then
we have
$\uu(\zz) \leqsw \uu(\zz_{-})$ by Lemma \ref{BumpLemPict} (1). 
Hence
$\zz \leqsw \zz_{-} \lsw \uu (\zz_{-})^{\tp} \leqsw \uu(\zz)^{\tp}$.
This proves Part (1). Part (2) follows immediately from Part (1).
$\hfill\square$

By \eqref{z0w1zkwk} and Part (2) of the lemma above,  there exists a unique cocorner 
$\ww_B$ of $\zeta$ such that 
\begin{equation}
\label{z-wBz+}
\zz_- \lsw \ww_B \lsw \zz_+.
\end{equation}
We will show that $\vv_B:= (\ww_B)^{\tp}$ is a balanced corner of $\LL$.

\begin{lemma}
The point $\vv_B$ is a removable corner of $\LL$ and satisfies
 $$
 \uu (\zz_-) \lsw \vv_B \lsw \uu (\zz_+).
 $$
\end{lemma}
\noindent
{\it Proof.} 
We give a proof of the second inequality.
The proof of the first inequality is similar and the first assertion follows from 
the second inequality and Lemma \ref{BumpLemPict} (6).   
Suppose that $\uu (\zz_+)$ is not an extreme cocorner of $(\mu/\zeta)^{\tp}$.
Then, we have $\uu (\zz_+)^{\tp} \in \ICC(\lz)$.
Since
$\zz_-$ is the predecessor of $\zz_+$ in $\bICC(\lz)$, 
$\zz_+ \in C_+$ implies $\uu (\zz_+)^{\tp} \leqsw \zz_-$. 
Hence we have $\vv_B \lsw \uu (\zz_+)$ by \eqref{z-wBz+}.
Next, suppose that $\uu (\zz_+)$ is an extreme cocorner of $(\mu/\zeta)^{\tp}$.
Since 
the inequality is obvious if $\uu (\zz_+) = \uu_{++}$, 
 it suffices to show that $\uu (\zz_+) \ne \uu_{--}$.
Suppose that $\uu (\zz_+) = \uu_{--}$ to the contrary.
By Lemma \ref{u--<uz++}, we have $\zz_+ \ne \zz_{++}$.
Hence $\zz_+ \lsw (\uu_{--})^{\tp} = \uu (\zz_+)^{\tp}$.
This contradicts to the fact that $\zz_+ \in C_+$.
$\hfill\square$

Define $\ww$ by $\displaystyle{  \vv_B\, \frac{}{\,\,\,{ \mathsf \Lambda}\,\,\,}\, \ww}$. 
Then 
we have $\zz_- \lsw \ww \lsw \zz_+$ 
by Lemma \ref{BumpLemPict} (5), (6) and the lemma above.
Since $\ww_B \in \IC (\lz)$ is characterized by \eqref{z-wBz+}, this proves
$\ww= \ww_B$. Thus, we get 
$\displaystyle{  \vv_B\, \frac{}{\,\,\,{ \mathsf \Lambda}\,\,\,}\, (\vv_B)^{\tp}}$, 
which completes the proof of Lemma \ref{exactness}.

\bigskip
\par\noindent
{\it Note.} After this paper was submitted to Journal of Algebra, another description of 
$m^{\mu}_{\lambda,(n-m,1^m)}$ was obtained by J. Blasiak: arXiv:1209.2018.



\end{document}